\documentclass[journal]{IEEEtran}
\usepackage{graphicx}
\usepackage{amssymb}
\usepackage{mathtools} 
\usepackage[margin=2.54cm]{geometry}
\usepackage{amsfonts} 
\usepackage[super]{nth} 
\usepackage{xcolor}
\usepackage{booktabs} 
\usepackage{multirow} 
\usepackage{subcaption} 
\usepackage{tabularx} 
\usepackage{float} 
\usepackage{xr}
\usepackage{bigstrut}
\usepackage{algorithm}
\usepackage{algorithmicx}
\usepackage{algpseudocode}
\usepackage{bbold}%

\ifCLASSINFOpdf
\else
\fi

\begin{document}
%
\title{Stochastic Process Optimization of an Integrated Biorefinery}
%
%
%

\author{Dahui~Liu$^1$, Sandra~Eksioglu$^1$,  Mohammad~Roni$^2$
\thanks{$^1$University of Arkansas, Fayetteville, AR, USA.}
\thanks{$^2$Idaho National Laboratory, Idaho Falls, ID, USA.}
\thanks{Manuscript received May 1, 2021.}}

%
%

\markboth{IEEE Transactions on Automation Science and Engineering, May~2021}%
{Shell \MakeLowercase{\textit{et al.}}: Bare Demo of IEEEtran.cls for IEEE Journals}
%



\maketitle

\begin{abstract}
Planning of biorefinery operations is complicated by the stochastic nature of physical and chemical characteristics of biomass feedstock, such as, moisture level and carbohydrate content.  Biomass characteristics affect the performance of the equipment which feed the reactor and the efficiency of the conversion process in a biorefinery. 
We propose a stochastic optimization model to identify a blend of feedstocks, inventory levels, and operating conditions of equipment to ensure a continuous flowing of biomass to the reactor while meeting the requirements of the biochemical conversion process. We propose a sample average approximation (SAA) of the model, and develop an efficient algorithm to solve the SAA model.  
A feedstock preprocessing process consists of two-stage grinding and pelleting is used to develop a case study. 
Extensive numerical analysis are conducted which lead to a number of observations. Our main observation is that sequencing  bales based on moisture level and carbohydrate content leads to robust solutions that improve processing time and processing rate of the reactor. We provide a number of managerial insights that facilitate the implementation of the model proposed.  
\end{abstract}

{\bf Note for practitioners:} This paper is motivated by the challenges faced in the bioenergy industry. The focus of this paper is on plants which use the biochemical conversion process to generate liquid fuels. It has been observed that variations in biomass characteristics, such as moisture  content, cause variations in feeding of the system  which lead to under-utilization of equipment.  A requirement of biochemical conversion process is to maintain the carbohydrate content of biomass processed by the reactor, larger than 59.1\%.  We propose a model that identifies the inventory levels and operating conditions of equipment to ensure a continuous flowing of biomass to the reactor. The goal is to improve equipment utilization while satisfying the requirements of the conversion process. The model is tested using real-life data. We found out that by sequencing bales based on moisture level and carbohydrate content, a plant can reduce variability in the system leading to higher processing rates of the reactor.    


%
\IEEEpeerreviewmaketitle

\section{Introduction} \label{introduction}
When the Renewable Fuel Standard (RFS) passed in 2007, it set annual target production amounts for cellulosic biofuel  which culminated to 10.5 Billion Gallons per Year (bgy) by 2020 \cite{RFS2020}. These annual targets for cellulosic biofuel were not met, what led the Environmental Protection Agency (EPA) to reduce the mandates. A 2016 study by US Department of Energy (DOE) indicates that US has the necessary resources (biomass) to produce up to 50 bgy of biofuels  by 2040 without jeopardizing the supply of food and animal feed \cite{langholtz20162016}. However, the cost of generating cellulosic biofuels is not competitive to that of fossil fuels. 

Based on \cite{RFS2020}, the shortfall in meeting the RFS is due to ``lack of private investment, logistical challenges, technology setbacks, and uneven support from the federal government."  A recent report from US DOE finds ``bulk solids handling and material flows through the system" to be a critical component to achieve the design throughput of the conversion processes \cite{DOE2016}. The flow of materials is impacted by the  variability of biomass characteristics which is observed when bales of different types of feedstock, with different moisture level and carbohydrate content, are processed by the same equipment (i.e., grinder). These variations negatively affect the integration of biomass feeding and conversion processes and cause low/unreliable on-stream time of equipment and low utilization efficiency of the reactor. Indeed, it was reported by Bell \cite{bell2005challenges} that inefficiencies of handling solid materials are the reason why biorefineries typically operate at 50\% of their capacity in the first year. Commercial scale cellulosic ethanol producers, such as, Beta Renewables (the 1st cellulosic ethanol plant in the world), Abengoa, DowDuPont, etc. have suspended their operations or filed for bankruptcy because cellulosic biofuel cannot compete with the low price of petroleum fuel. 

The objective of the proposed research is to develop operational strategies which streamline processes within a biorefinery to ensure a continuous flow of biomass to the reactor, while meeting the requirements of biochemical conversion process. This will lead to shorter processing times, higher utilization of the reactor, and lower processing costs. To accomplish this, we propose a mathematical model that identifies the infeed rate of biomass to the system, processing speed of equipment and inventory level to ensure a continuous feeding of the reactor. It is assumed that different biomass feedstocks are used in the system, each with different moisture and carbohydrate contents. Moisture level impacts the infeed rate and processing speed of equipment since biomass of high moisture level is processed slow to avoid equipment clogging. Carbohydrate content depends on the feedstock type. The biochemical conversion process requires that the carbohydrate content of biomass should be at least 59.1\%. Since, both, moisture and carbohydrate content are random variables, we propose a stochastic optimization model.  

Based on our numerical analysis, sequencing  bales based on their moisture and carbohydrate contents  reduces variations of the characteristics of biomass flowing in the system, which could positively impact system's performance. Integrating a bale sequencing model within the proposed stochastic optimization model would make the model challenging to solve. Therefore, we evaluate different sequences numerically. Via our numerical analysis we identify characteristics of good  sequences and develop a few bale sequencing rules that are easy to implement. We also provide a number of managerial insights to facilitate large-scale implementations of the models proposed.  

We develop a case study using data collected at the Process Development Unit (PDU) of Idaho National Laboratory (INL). PDU is a full-size, fully-integrated feedstock preprocessing system. A number of INL's industry partners use the PDU to test different grinding, pelleting, mechanical and chemical separation options, etc. during design and scale-up of a biorefinery. Figure \ref{flow chart} presents the flowchart of processes of PDU. Current processes include drag chain conveyor 1 (DC(1)) to pelleting mill (P). The proposed model also considers a separation process, which separates biomass based on particle size. Large particles are processed further. The proposed model identifies a blend of biomass to meet biochemical conversion specifications as designed by NREL \cite{davis2013process}. These specifications are a total structural carbohydrate of 59.1\%, and a moisture content of 20\%.  

In summary, these are the main contributions of this research: ($i$) a stochastic optimization model is developed that leads to robust process control solutions that meet biochemical process requirements and maintain a continuous flow of biomass to the reactor; ($ii$) a real-life case study is developed using historical data from INL's PDU, and other reliable data sources; ($iii$) a number of managerial insights are identified that have the potential to positively impact the performance of biomass processing systems.   

The reminder of the paper is organized as follows. In Section \ref{literature review} we briefly review the literature. The proposed stochastic optimization model is introduced in Section \ref{model}. In Section \ref{approach} we present the algorithms developed. In Section \ref{experiment} we present our experimental setup. The results of our experiment and observations we make are summarized in Section \ref{result}. In Section \ref{insight} we provide managerial insights about the implementation of the proposed model.

\section{Literature Review} \label{literature review}
The proposed work aligns with two main streams of literature within the broad area of  supply chain optimization. First, the proposed model is an extension of the  blending problem which identifies minimum cost biomass blends to meet process requirements of conversion process. Second, the proposed work contributes to the existing supply chain reliability literature.  Below, we provide details of the literature and highlight our contributions.

\subsection{Biomass Blending}
The \emph{blending model} is a well-known model in operations research \cite{winston2004operations}. It identifies a blend of resources to meet demand (or some system requirement) to minimize costs (or to maximize profits). The model is used in a number of application in different industries. For example, the model is used in the agricultural sector to adjust the quality (grade) of grains sold in the grain market \cite{hill1990grain}. The model is used in the energy sector to identify coal blends to meet meet emission criteria \cite{sami2001co, shi2015impact, Huang2009, Wang2018}. Most recently, a number of researchers are investigating the opportunities that exist to reduce the cost of producing biofuels by using biomass blends in biorefineries \cite{roni2018optimal,sharma2020assessment,eksioglu2020}.

The research on biomass blending is motivated by the fact that biomass feedstocks  differ based on their physical characteristics, chemical content, and selling price. There is an opportunity to reduce the cost of producing bioenergy by using a blend of feedstocks that meet process requirements at the minimum cost. However, using biomass blends in a biorefinery is challenging for the following reasons  ($i$) processes should be re-evaluated and re-adjusted to handle a blend of products with different physical properties, which impact biomass flowability in the system; ($ii$) the efficiency of the conversion process is affected by the chemical composition of the biomass blend; ($iii$) business relationships should be built and maintained with a number of suppliers to obtain feedstocks; etc. In response to these challenges,  three main streams of research have been developed that focus on evaluating conversion process \cite{crawford2015evaluating, ray2017biomass, wolfrum2017effect}, evaluating the economic feasibility and risks within the biorefinery \cite{thompson2017optimizing, ou2018understanding, lan2020impacts, lan2021techno}, and evaluating the supply chain costs and risks of using biomass blends  \cite{jacobson2014biomass,roni2018optimal, narani2019simultaneous, roni2019,  chen2020ensemble, lan2020life}.

The proposed research contributes to the literature focused on evaluating the economic feasibility of biomass blend within a biorefinery. Different from the literature, the proposed work demonstrates the impact of biomass blending on flow, reactor's utilization, inventory level, etc. We propose strategies to sequence bales which ensure a continuous flow of biomass to the reactor.   

\subsection{Systems Reliability}
A recent work by \cite{DavidCoit_reliability_opt2018} provides a review of the reliability literature and classifies existing models into ($i$) redundancy allocation, ($ii$) reliability allocation, ($iii$) reliability-redundancy allocation, and ($iv$) assignment and sequencing. Our proposed model contributes to redundancy allocation and sequencing streams of research since it identifies the amount of inventory to ensure continuous feeding of the reactor. In the numerical analysis, we investigate the impact of bale sequencing in the performance of the system. 

The model proposed falls under ($i$) redundancy allocation and ($iv$) assignment and sequencing streams of research for the following reasons. Our model identifies buffer size and location, that add redundancies to the system. This inventory contributes to maintaining a continuous flow of biomass to the reactor. We have  evaluated numerically the impact of bale sequencing based on moisture level and carbohydrate content to the performance of the system.  Other researchers use extensions of the redundancy allocation model to design  processes in a bioenergy plant \cite{andiappan2016synthesis, andiappan2017systematic, andiappan2019integrated,  gulcan2021optimization, Chebouba2018}. Different from these works, our model considers the impact of bale sequencing, and biomass blending to processes in a bioenergy plant.

The model we propose is an extension of the work by \cite{kucuksayacigil2021optimal} that identifies a process control strategy to streamline  preprocessing of switchgrass in a biorefinery. Different from  \cite{kucuksayacigil2021optimal}, we propose a stochastic optimization model that considers the impacts of stochastic carbohydrate and moisture content to the performance of the system. Additionally, our proposed model considers a blend of biomass feedstocks.  
 
Finally, the model we propose is an extension of the flow shop  scheduling problem that identifies a sequence of (bales) jobs processed that minimizes makespan \cite{GJS76,graham1979optimization,zhang2019review, Mhasawade2017, Ni2020}. The order of operations is the same for all (bales) jobs, but the processing time of each (bale) job  in each operation differs based on the characteristics of the (bale) job. In our model, we group bales based on moisture level (low, medium, high) and carbohydrate content ($< 59.1\%$ or $\geq 59.1\%$). Since this grouping leads to only a few different types of bales, we evaluate sequences numerically.   

\section{Model Formulation} \label{model}
\subsection{Problem Description}
\begin{figure*}[htb]
    \centering
    \includegraphics[width=1.2\columnwidth]{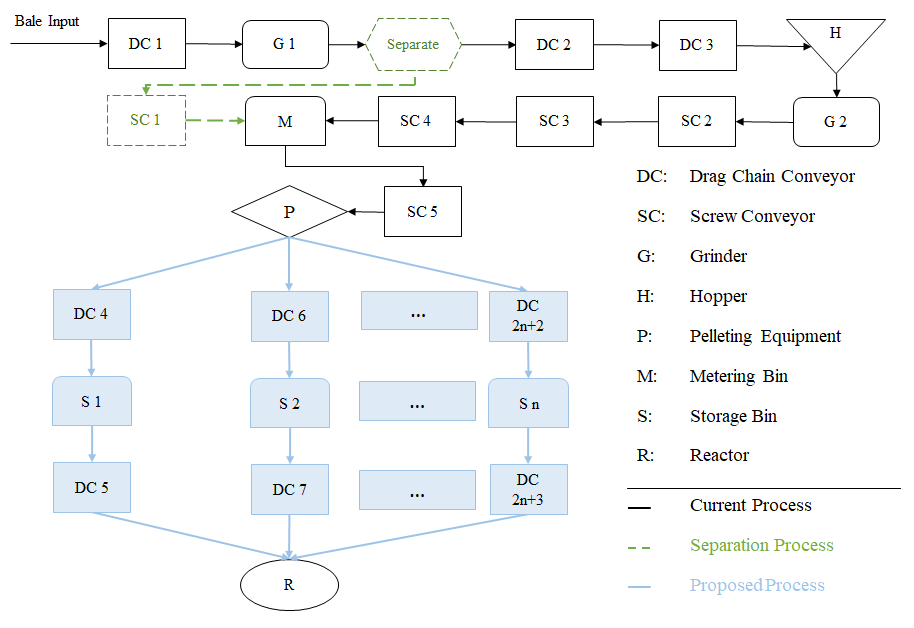}
    \caption{Proposed Flowchart of Processes at PDU.}
    \label{flow chart}
\end{figure*}

Figure \ref{flow chart} presents a flowchart of processes at the PDU. The current processes include DC(1), which feeds biomass bales to the system, to P where biomass is pelleted. DC(1) moves biomass to Grinder 1 (G(1)). A separation process after G(1) separates biomass based on particle size. Large particles are moved via DC(2) and DC(3) to G(2) for further processing, and from G(2) processed biomass is moved to metering bin (M). Small particles are moved via SC(1) to M. Via SC(5), biomass is moved from M to pelleting mill (P). 

The proposed flowchart includes additional storage (S(1)  to S($n$)) and conveyors (DC(4) to DC(2$n$ +3)). Each additional storage unit inventories pellets of a particular biomass feedstock. Each feedstock has specific carbohydrate and ash contents. A mix of pellets, which meet biochemical conversion process requirements, is fed to the reactor (R) in every time period.  

\subsection{A Mixed-Integer Nonlinear Model}
Let us assume that $n_{sb}$ bales of feedstock $b\in\mathcal{B}$, with moisture level $s \in \mathcal{S}$, weighting $m_{sb}$ dry tons (dt) each, are to be processed at the PDU. Let $t (\in \mathcal{T})$ represent a time period within the planning horizon. The objective is to identify how much biomass to process in every time period and how much biomass to carry in the inventory in order to minimize the time it takes to process every bale. Let $\tilde{Z}_t$ be a binary variable that takes the value 1 if the reactor is up and running in period $t$, and takes the value 0 otherwise. The objective of our model is
\begin{align}
\min &  \sum_{t \in \mathcal{T}} \tilde{Z}_{t}. \label{objective_0}
\end{align} 

Biomass processing is impacted by the capacity of equipment used. Let $u_{isb}$ represent the capacity of equipment $i (\in \mathcal{I})$ in terms of dt of biomass processed per time period. The capacity depends on the type of biomass processed and its moisture level. For example, the PDU operator lowers the rotational speed of grinders when processing highly moist biomass in order to avoid equipment clogging. Thus, the processing capacity of an equipment is higher for dry biomass then for moist biomass. Additionally, the capacity of storage bins ($\mathcal{I^{M}} \cup \mathcal{I^{B}}$) is expressed in terms of mass ($m_i$) and volume ($v_i$), because the density of processed biomass ($d_{isb}$) impacts the space it takes in storage. Biomass feedstocks differ based on density, which impacts the amount of biomass that can be stored in the same storage space, thus, equipment capacity.  

Let $X_{isbt}$ represent the amount of biomass flowing from equipment $i$, and let $M_{isbt}$ represent the amount of biomass in inventory. \eqref{eq:outflowmax0} to \eqref{eq:volumecapacityexpanded0}  represent the capacity constraints of processing, transportation and storage equipment. 
\begin{align}
X_{1sbt} \leq {u}_{1sb} {Z}_{sbt}, & \forall s \in \mathcal{S}, b \in \mathcal{B}, t \in \mathcal{T}, \label{eq:outflowmax0} \\
X_{isbt}  \leq {u}_{isb}, & \forall i \in \mathcal{I},  s \in \mathcal{S},  b \in \mathcal{B},  t \in \mathcal{T}, \label{eq:outflowmaxG0} \\
\sum_{s \in \mathcal{S}} \sum_{b \in \mathcal{B}} M_{isbt} \leq {m}_{i}, & \forall i \in \mathcal{I^{M}} \cup \mathcal{I^{B}},  t \in \mathcal{T}, \label{eq:masscapacityexpanded0} \\
\sum_{s \in \mathcal{S}} \sum_{b \in \mathcal{B}} \frac{M_{isbt}}{{d}_{isb}}  \leq {v}_{i}, & \forall i \in \mathcal{I^{M}} \cup \mathcal{I^{B}}, t \in \mathcal{T}. \label{eq:volumecapacityexpanded0} 
\end{align}
In constraint \eqref{eq:outflowmax0}, ${Z}_{sbt}$ is binary decision variable which takes the value 1 if biomass feedstock $b$ is processed in period $t$, and takes the value 0 otherwise. The value of these variables is determined by constraint \eqref{eq:onlySBcanprocess0}.  When ${Z}_{sbt} = 0,$ the feeding of the system ends, what  eventually ends reactor's processing. 
\begin{align}
\sum_{s' \in \mathcal{S}} \sum_{b' \in \mathcal{B}} \sum_{t'=t}^{t+p_{sb}-1} {Z}_{s'b't'} - \sum_{t'=t}^{t+p_{sb}-1} {Z}_{sbt'} & \leq M (1 - y_{sbt}), \nonumber \\ \forall s \in \mathcal{S},  b \in \mathcal{B},  t \in \mathcal{T}. \label{eq:onlySBcanprocess0} 
\end{align}

In  constraint  \eqref{eq:onlySBcanprocess0}, $y_{sbt}$ is a parameter which takes the value 1 if a bale with characteristics $b$ and $s$ is  processed in period $t$, and takes the value zero otherwise. The values of this parameter represent the sequence of bales fed to the system. In this model we consider this to be a problem input. Nevertheless, mathematical models can be used to sequence the bales. 

The following constraints mimic how the system operates. Constraint \eqref{eq:masspertimeperiod0} determines the amount of biomass fed to the system in period $t$. This amount depends on the infeed rate of the system (in meters per time unit), $V_{sbt}$, and the amount of biomass in (each meter of) a bale, $c_{sb}$. To estimate $c_{sb} (= w\times h \times d_{sb})$ we use the dimensions of a bale ($w, h$) and its density ($d_{sb}$). Constraint \eqref{eq:allprocessed0} ensures that the amount of biomass fed to the system during the planning horizon does not surpass the amount available.  Constraints \eqref{eq:movedistancegreaterl0} and \eqref{eq:movedistancelessl0} present the relationship between the infeed rate and the time it takes to process a bale of length $l$. Constraint \eqref{eq:become0keep1} determines the values of $\tilde{Z}_t,$ and helps us determine the minimum processing time of the reactor.
\begin{equation}
    X_{1sbt}  \leq c_{sb} V_{sbt}, \forall s \in \mathcal{S},  b \in \mathcal{B},  t \in \mathcal{T},  \label{eq:masspertimeperiod0}
\end{equation}
\begin{equation}
   \sum_{t \in \mathcal{T}} X_{1sbt}  \leq m_{sb} n_{sb}, \forall s \in \mathcal{S},  b \in \mathcal{B}, \label{eq:allprocessed0} 
\end{equation}

\begin{equation}
   \sum_{t'=t}^{t+p_{sb}-1} V_{sbt'}  + {M} ( 1-y_{sbt} ) \geq l, \forall s \in \mathcal{S},  b \in \mathcal{B},  t \in \mathcal{T}, \label{eq:movedistancegreaterl0}
\end{equation}
\begin{equation}
  \sum_{t'=t}^{t+p_{sb}-1} V_{sbt'} - {M} ( 1-y_{sbt} ) \leq l, \forall s \in \mathcal{S},  b \in \mathcal{B},  t \in \mathcal{T}, \label{eq:movedistancelessl0}
\end{equation}

\begin{equation}
\tilde{Z}_{t} \leq \tilde{Z}_{t-1} \forall t \in \mathcal{T}\backslash\{1\}. \label{eq:become0keep1}
\end{equation}

Constraints \eqref{eq:flowgrinders0} to \eqref{eq:flowdrags0} represent the flow balance constraints. In  \eqref{eq:flowgrinders0}, the amount of flow from  grinder $i$ is  reduced based on dry matter loss, $\mu_i$ (\%). Constraints \eqref{eq:flowdcc50} and \eqref{eq:flowsc60} model the separation process. Here, $\theta_{sb}$ is the bypass ratio of this system. 

\begin{equation}
X_{isbt}  = \sum_{p \in \mathcal{I}_{i}} (1 - \mu_{i}) X_{psbt},    \forall i \in \mathcal{I^{G}},  s \in \mathcal{S},  b \in \mathcal{B},  t \in \mathcal{T}, \label{eq:flowgrinders0} 
\end{equation}

\begin{align}
X_{ksbt}  = \sum_{p \in \mathcal{I}_{i}} \left( 1 - \theta_{sb} \right) X_{psbt},  
k = \text{DC(2)}, \forall s \in \mathcal{S},  b \in \mathcal{B},  \nonumber\\  t \in \mathcal{T}, \label{eq:flowdcc50}
\end{align}

\begin{equation}
X_{ksbt}  = \sum_{p \in \mathcal{I}_{i}} \theta_{sb}  X_{psbt},  k = \text{SC(1)}, \forall s \in \mathcal{S},  b \in \mathcal{B},  t \in \mathcal{T}, \label{eq:flowsc60} 
\end{equation}

\begin{align}
X_{isbt}  = \sum_{p \in \mathcal{I}_{i}} X_{psbt}   \forall i \in \mathcal{I}\backslash\{ \text{DC(2), SC(1)}, \mathcal{I^M}, \mathcal{I^{B}}, \mathcal{I^G} \}, \nonumber\\ s \in \mathcal{S},  b \in \mathcal{B},  t \in \mathcal{T}.\label{eq:flowdrags0}
\end{align}

Constraint \eqref{eq:invall0} presents the flow conservation constraint in storage bins. Constraints \eqref{eq:initialzero0} and \eqref{eq:allOut0} show that in initial and ending inventory are zero. 
\begin{align}
M_{isbt}  = M_{isbt-1} + \sum_{p \in \mathcal{I}_{i}} X_{psbt} - X_{isbt} & \forall i \in \mathcal{I^{M} \cup I^{B}}, \nonumber\\ s \in \mathcal{S},  b \in \mathcal{B},  t \in \mathcal{T}, \label{eq:invall0} \\
M_{isb0}  = 0, \forall i \in \mathcal{I^{M} \cup I^{B}},  s \in \mathcal{S},  b \in \mathcal{B}, \label{eq:initialzero0}\\
M_{isbT}  = 0,  \forall i \in \mathcal{I^{M} \cup I^{B}},  s \in \mathcal{S}, b \in \mathcal{B}. \label{eq:allOut0}
\end{align}

Constraints \eqref{eq:100utilization0} to \eqref{eq:aveutilization0} model reactor's reliability. These constraints ensure that the feeding of the reactor is within a certain  bound of its maximum capacity. Let $U$ represent the maximum feeding rate of the reactor. Constraints \eqref{eq:100utilization0} and \eqref{eq:minutilization0} ensure that the feeding of the reactor is within $[\underline{q}U, U]$ (for $0 \leq \underline{q} < 1$) in every time period. Additionally, \eqref{eq:aveutilization0}  ensures that, on the average, the feeding of the reactor is at least $qU$ for $q \in (\underline{q}, 1]$. 
We multiply the right-hand-side of \eqref{eq:minutilization0} and \eqref{eq:aveutilization0} by $\tilde{Z}_t$ to estimate the amount of biomass processed during the time that the reactor is up and running (i.e., $\tilde{Z}_t =1$).  
\begin{align}
\sum_{i \in \mathcal{I}_r} \sum_{s \in \mathcal{S}} \sum_{b \in \mathcal{B}} X_{isbt}, & \leq {U}, & \forall t \in \mathcal{T}, \label{eq:100utilization0} \\
\sum_{i \in \mathcal{I}_r} \sum_{s \in \mathcal{S}} \sum_{b \in \mathcal{B}} X_{isbt} & \geq \underline{q} {U}\tilde{Z}_t, & \forall t \in \mathcal{T}, \label{eq:minutilization0} \\
\sum_{i \in \mathcal{I}_r} \sum_{s \in \mathcal{S}} \sum_{b \in \mathcal{B}} \sum_{t \in \mathcal{T}} X_{isbt} & \geq {q} U \sum_{t \in \mathcal{T}} \tilde{Z}_t. \label{eq:aveutilization0}
\end{align}

Let $\mathcal{X}$ represent the feasible region of the optimization model described via constraints \eqref{eq:outflowmax0} to \eqref{eq:aveutilization0}, and ${\bf x} \in \mathcal{X}$ a solution in this region. Let $f({\bf x})$ represent the objective \eqref{objective_0} as a function of ${\bf x}$. The following is a succinct formulation of the optimization model described in this Section
\[\min: f({\bf x}) \mbox{ s.t. } {\bf x}\in \mathcal{X}.\]
Let us refer to this as model formulation ($P$).

\subsection{A Linear Reformulation of Model ($P$)} \label{LinearSAA}
In  ($P$), the maximum feeding rate of the reactor, $U$, is a decision variable. We do not set this  variable to a predetermined value in order to allow the model identify the maximum feeding rate of the reactor for given the available storage space, the quality of biomass and the processing speed of equipment.  

In constraints \eqref{eq:minutilization0} and \eqref{eq:aveutilization0} of ($P$) the terms ${U}\tilde{Z}_t$ are bilinear, thus,  ($P$) is a  mixed-integer nonlinear program (MINLP). We use the McCormick relaxation \cite{mccormick1976computability} of these bilinear terms. The corresponding relaxation is  a mixed integer linear program (MILP). Furthermore, since $\tilde{Z}_t$ are binary variables, this relaxation is tight \cite{Adams1993}. 

Let ($\bar{P}$) be the new MILP. In ($\bar{P}$), $W_t = {U}\tilde{Z}_t$. Thus,  the right-hand-side of \eqref{eq:minutilization0} and \eqref{eq:aveutilization0} are updated to $\underline{q}W_t$ and $q\sum_{t\in\mathcal{T}} W_t$, respectively. Additionally, the following constraints are added to the constraint set.
\begin{align}
W_t & \geq \underline{U}\tilde{Z}_t, \hspace{1.1in}\forall t\in \mathcal{T},\label{mc1}\\
 W_t & \leq \overline{U}\tilde{Z}_t, \hspace{1.1in}\forall t\in \mathcal{T},\label{mc2}\\
W_t & \geq \overline{U}\tilde{Z}_t + U - \overline{U}, \hspace{0.5in}\forall t\in \mathcal{T},\label{mc3}\\
W_t &\geq \underline{U}\tilde{Z}_t + U - \underline{U},  \hspace{0.5in}\forall t\in \mathcal{T}.\label{mc4}
\end{align}
In constraints \eqref{mc1} to \eqref{mc4}, $\underline{U}, \overline{U}$ are fixed and known lower and upper bounds of $U.$  Let $\bar{\mathcal{X}}$ represent the feasible region of ($\bar{P}$).  

\subsection{A Chance Constraints Model} \label{SAA}
Blending of biomass feedstock with different physical and chemical characteristics provides an opportunity to meet the requirements of the conversion platform at minimum cost. For example, a  report by the National Renewable Energy Laboratory (NREL) \cite{davis2013process} sets the design specifications for the biochemical conversion. These specification recommend that the total structural carbohydrate of the feedstock used should be at least 59.1\%.  

\begin{figure*}[h!t]
    \centering
    \includegraphics[width=0.8\textwidth]{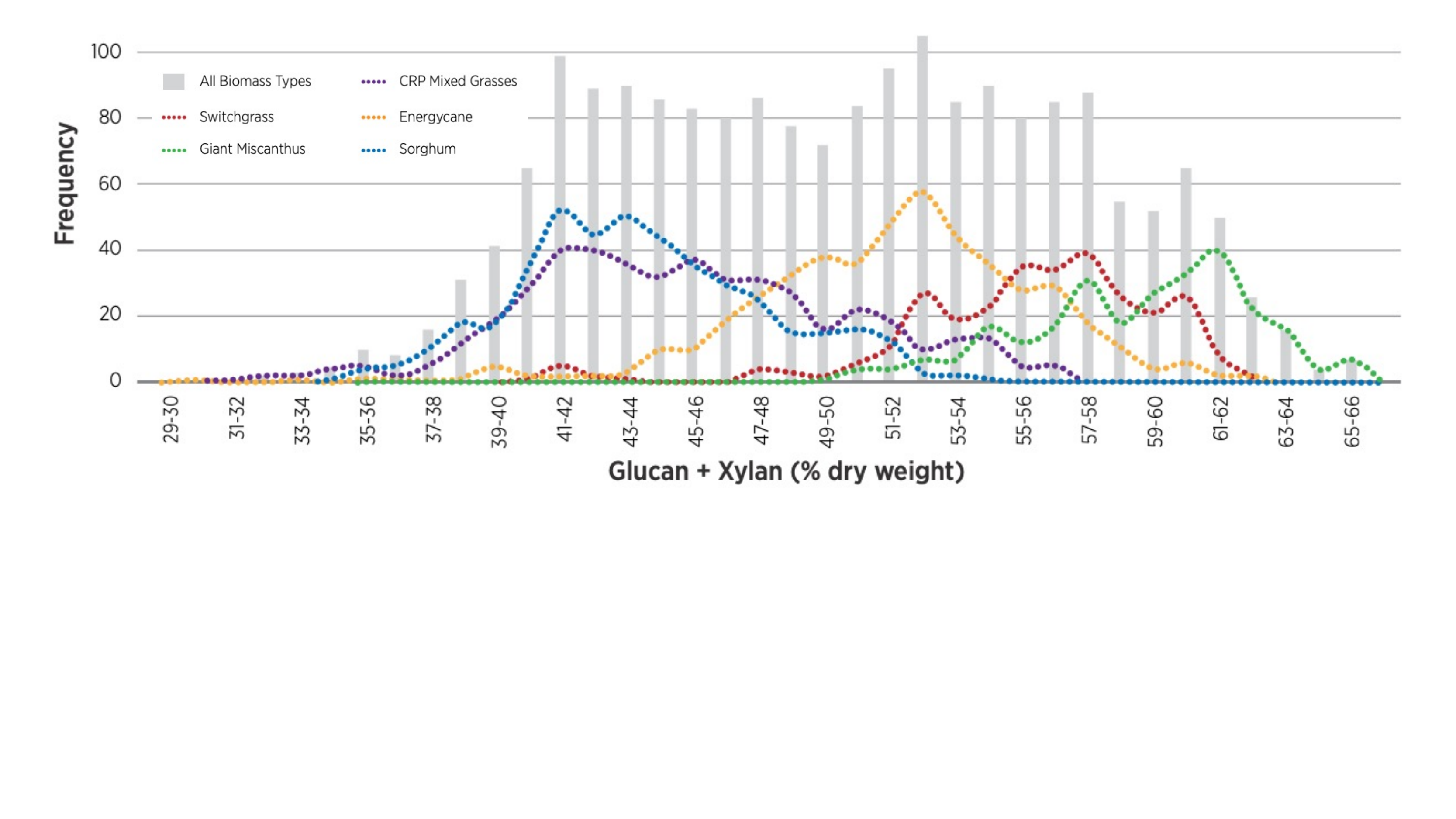}\vspace{-1.2in}
    \caption{Distribution of Carbohydrates (\% glucan plus \% xylan) for Biomass Feedstocks \cite{Owens2016RegionalFP}.} \label{fig-carb}
\end{figure*}

A total of 2,000 samples of biomass was collected and analyzed during 2008 and 2014 by \cite{Owens2016RegionalFP}. A summary of the carbohydrate contents of these samples is presented in Figure \ref{fig-carb}. Based on this data, switchgrass and mischanthus have higher carbohydrate contents as compared to other feedstocks analyzed. The average carbohydrate contents for these feedstocks are close to the 59.1\% threshold, which means that about 50\% of the samples collected do not meet the threshold.    Although biorefineries would prefer to use biomass feedstocks with carbohydrate contents higher then this threshold all the time, achieving this may be impractical and costly. Thus, in our model setting we assume that a biorefinery strives to meet the carbohydrate content threshold most (i.e., 80\% or 90\%) of the time. To model this, we use chance constraints. 

Let $f_b$ be a stochastic parameter which represents the carbohydrate content of biomass feedstock $b$. Let $f^*$ represent the threshold  level that the biorefinery strives to achieve, and let $\gamma$ represent the  risk level associated with achieving this threshold. Constraint \eqref{eq:blendingf} is the chance constraint which indicates that, during the planning horizon $\mathcal{T}$, the carbohydrate threshold requirements are met  at least $(1-\gamma)\%$ of the time.
\begin{align}
Pr \left( \sum_{i \in \mathcal{I}_{r}} \sum_{s \in \mathcal{S}} \sum_{b \in \mathcal{B}} \sum_{t \in \mathcal{T}} (f_{b} - f^*) X_{isbt} \geq 0 \right ) & \geq 1 - \gamma. \label{eq:blendingf} 
\end{align}
Let $H ({\bf x},f) = \sum_{i \in \mathcal{I}_{r}} \sum_{s \in \mathcal{S}} \sum_{b \in \mathcal{B}} \sum_{t \in \mathcal{T}} (f^*-f_{b}) X_{isbt}.$ The following is a succinct formulation of the chance constraint optimization model, which we refer to as model ($Q$).  
\[\min_{\bf x\in \mathcal{X}}: f({\bf x}) \mbox{ s.t. } Pr\left(H ({\bf x},f) \geq 0\right) \leq \gamma.\]
\subsection{A Sample Average Approximation Model}
Chance constraint models, such as ($Q$), are difficult to solve for two main reasons \cite{Birge_Louveaux_1997}. First, calculating the probability of meeting a certain constraint is often challenging due to multidimensional integration. Second, the feasibility region defined by the chance constraints may not be convex \cite{kim2015guide}. There are two main approaches solve chance constraints models. Both approaches approximate the probabilistic constraints. One approach discretizes the probability distribution function of the random parameters, and solves the corresponding combinatorial problem \cite{Dentcheva_Prekopa_Ruszczynski_2000, Luedtke_Ahmed_2008}. The other approach develops convex approximations of the probabilistic constraint \cite{Nemirovski_Shapiro_2006}. We follow the first approach to solve ($Q$).

In model ($Q$), the probability density function of $f_b$ does not follow a known distribution (see Figure \ref{fig-carb}), thus, calculating $Pr\left(H ({\bf x},f) \geq 0\right)$ is challenging. We use the Sample Average Approximation (SAA) method to approximate the probability density function of $f_b$ via an empirical density function.  Let $\hat{H}_N({\bf x},f) = N^{-1}\sum_{n=1}^N\Delta(H_n({\bf x}, f_{n}))$ denote an empirical measure of the probability distribution of $H({\bf x},f).$ Here, $f_{n}=\{f_{bn}| b\in \mathcal{B}\}$ (for $n=1,\ldots,N$) are $N$ iid samples of $f_b$; $\Delta()$ is a measure of the probability mass function value; and $N^{-1}$ is the probability assigned to  observations of $f_b$.  

We can approximate $H({\bf x},f)$ using the empirical measure $\hat{H}_N({\bf x},f)$ as follows.
\begin{align}\label{eq-estimate-H}
\bar{H}_N({\bf x}) = \mathbb{E}\bigg[\mathbb{1}_{(0,\infty)} \big(\hat{H}_n({\bf x}, f_{n})\big)\bigg] = \nonumber \\ \frac{1}{N}\sum_{n=1}^N\mathbb{1}_{(0,\infty)} (\hat{H}_n({\bf x}, f_{n})).
\end{align}
Where, $\mathbb{1}_{(0,\infty)}(t) := \bigg{\{} \begin{array}{c}  1 \mbox{ if } t > 0\\ 0 \mbox{ if } t \leq 0 \end{array}.$ $\bar{H}_N({\bf x})$ returns what proportion of the time  $\hat{H} ({\bf x},f) \geq 0.$  The resulting SAA approximation of ($Q$), referred to as model ($\bar{Q}$), is presented below. 
\[\min_{\bf x\in \mathcal{X}}: f({\bf x}) \mbox{ s.t. } \bar{H}_N({\bf x}) \leq \hat{\gamma}.\]

Note that $\hat{\gamma}$ is the risk level  of ($\bar{Q}$). Based on \cite{Luedtke_Ahmed_2008}, if $\hat{\gamma} < \gamma$, then, the probability that a feasible solution of model ($\bar{Q}$) is feasible to the true model ($Q$) approaches 1 as $N\rightarrow \infty.$

\section{Solution Approach} \label{approach}
\subsection{An Algorithm to Solve Model ($\bar{Q}$)}
The constraint set of ($\bar{Q}$) includes an indicator function. Since commercial software cannot handle such a function, we reformulate these constraints by introducing the following continuous variables, $\beta^{-}$ and $\beta^{+}.$ These variables quantify the violation of constraint $\bar{H}_N({\bf x}) \leq \hat{\gamma}$. In the following model we have approximated this constraint using linear functions. 
\begin{align}
\nu^*(\alpha) = \min_{\bf x\in \mathcal{X}}:  f({\bf x})  + \alpha \sum_{n =1}^N  \beta_{n}^{-} \label{objective_3}
\end{align} \vspace{-0.2in}
\begin{align}
\text{Subject to:}\hspace{1in}&&\nonumber \\
{H}({\bf x},f_n) + \beta_{n}^{+} - \beta_{n}^{-} & = 0 &  \forall n=1,\ldots, N, \label{Q_hat_1}\\
\beta_{n}^{-},\beta_{n}^{+} & \geq 0 &  \forall n=1,\ldots, N. \label{Q_hat_2}
\end{align}
We refer to this formulation as model ($\hat{Q}_{\alpha}$). Its objective minimizes the cost of violating the chance constraint. The value of the parameter $\alpha$ in the objective function determines the penalty for violating the chance constraint. When the value of $\alpha$ is very large, the objective sets the value of $\beta^{-}_n = 0$ for all $n=1,\ldots, N.$  Consequently, $H_n({\bf x},f_n) \leq 0$ for all $n=1,\ldots, N.$  When the value of $\alpha$ is very small, the objective sets the value of $\beta^{-}_n \geq 0$ which leads to solutions which violate the chance constraint. Thus, the goal is to find $0 <  \alpha \in \big[\alpha^{l}, \alpha^{u}\big]$ which minimizes $\nu^*(\alpha)$ while ensuring that the chance constraints are satisfied in at least $\lceil(1- \hat{\gamma})N\rceil$ of the samples generated.
\begin{align}
\min_{\alpha}: \nu^*(\alpha)\label{objective_4}
\end{align}\vspace{-0.2in}
\begin{align}
\text{Subject to:}\hspace{2in}&&\nonumber \\
0 <  \alpha \in \big[\alpha^{l}, \alpha^{u}\big],\\
\sum_{n=1}^N\mathbb{1}_{(0,\infty)} \left(\beta_n^{-}\right) \geq \lceil(1- \hat{\gamma})N\rceil.  
\end{align}
We refer to this formulation as model ($\hat{Q}$). The SAA Algorithm described below solves ($\hat{Q}$) using a binary search approach.
\begin{algorithm}\label{SAA-alg}
	\caption{SAA Algorithm}
	\textbf{Initialization:} $0 < \alpha^{l} \leq \alpha \leq \alpha^{u}$; $\epsilon = 10^{-4},$ $\phi = 10^{-4}$; {\bf Found} = False.
	\begin{algorithmic}[1]
		\While {({\bf Found} = False)}
	        \State $\alpha \leftarrow \frac{\alpha^l + \alpha^u}{2}$
		    \State Solve ($\hat{Q}$); Let $\hat{\beta}_{n}^{-}$ be the incumbent solution of $\beta_{n}^{-}$ 
		        \State $C \leftarrow 0$
		        \For {$n = 1,\ldots,N$}
		            \If {$\beta_{n}^{-} > 0$}
		                \State $C \leftarrow C + 1$
		            \EndIf
		        \EndFor
		        \If {$C \geq \hat{\gamma} N + \epsilon$}
		            \State $\alpha^l \leftarrow  \alpha$
		        \ElsIf {$C < \hat{\gamma} N - \epsilon$}
		        \State $\alpha^u \leftarrow  \alpha$
		        \EndIf
			\If {$|\alpha - \frac{\alpha^{l} + \alpha^{u}}{2} | \leq \phi$}
			    \State {\bf Found} = True 
			\EndIf
		\EndWhile
		\State Return $\alpha$ and the solution of ($\hat{Q}$) 
	\end{algorithmic} 
\end{algorithm} 

\subsection{An Algorithm to Solve Model ($Q$)}\label{Sec:UpperLowerBd}
Let us briefly describe the algorithms we use to generate probabilistic lower and upper bounds for the true chance constraint model ($Q$). We follow closely the approach proposed by \cite{Luedtke_Ahmed_2008}. 

Let $\nu^*$ and $\hat{\nu}_{N}$ represent the optimal objective function value of ($Q$) and  ($\hat{Q}$), respectively. Since ($\hat{Q}$) is a relaxation of ($Q$), and for $\hat{\gamma} > \gamma$,  $\hat{\nu}_{N}$ is a lower bound for $\nu^*$. Indeed, the probability that $\hat{\nu}_{N}$ is a lower bound to $\nu^*$ approaches 1 as $N \rightarrow \infty.$ If $\hat{\gamma} < \gamma$, then, the probability that $\hat{\nu}_{N}$ is an upper bound for $\nu^*$ approaches 1 as $N\rightarrow \infty.$ Work by \cite{Luedtke_Ahmed_2008} determines the sample size $N$, the number of replications $M$, and the following approaches to find high confidence upper/lower bounds for ($Q$).

\emph{Lower Bound Procedure}: 
We follow these steps: (1) Select $\hat{\gamma} > \gamma$, and $N\geq \frac{1}{2(\hat{\gamma} - \gamma)^2}\log{\big(\frac{1}{\delta}\big)}$. (2) Solve $M=10$ replications of ($\hat{Q}$). In each replication we use a different sample of size $N.$ Let $\hat{\nu}_{Nj}$ represent the corresponding objective function values for $j=1,\ldots, M.$ (3) Find the $\underline{\nu}= min_{j=1,\ldots,M} \hat{\nu}_{Nj}.$ Then, $\underline{\nu}$ represents the lower bound of ($Q$) with confidence $(1-\delta)$. 

\emph{Upper Bound Procedure}: We follow these steps: (1) Select $N$, $\delta$, $\hat{\gamma} < \gamma$, $M=10$. (2) Solve ($\hat{Q}$) for each of the $M$ replications, and let ${\bf x}_{jN}$ represent the solution found ($j=1,\ldots, M$). (3) Conduct an \emph{a posteriori} check to see if the chance constraint is satisfied, that is, $Pr\left(H ({\bf x}_{jN},f) \geq 0\right) \leq \gamma.$ In order to estimate this probability, we sample $N^{\prime}$ ($>>N$) iid values of the random parameter $f_b.$ We use equation \eqref{eq-estimate-H} to calculate $\bar{H}_{N^{\prime}}({\bf x}_{jN},f)$, which is an estimate of the  $Pr\left(H ({\bf x}_{jN},f) \geq 0\right)$. (4) Calculate $\mathcal{U}_{N^{\prime}}({\bf x}_{jN},f) = \bar{H}_{N^{\prime}}({\bf x}_{jN},f) + \Phi^{-1}(1 - \delta)\sqrt{\frac{\big(\bar{H}_{N^{\prime}}({\bf x}_{jN},f)\big)\big(1-\bar{H}_{N^{\prime}}({\bf x}_{jN},f)\big)}{N^{\prime}}},$ where $\Phi^{-1}()$ represents the inverse c.d.f. of $N(0,1)$. If $\mathcal{U}_{N^{\prime}}({\bf x}_{jN},f) \leq \gamma,$ then, ${\bf x}_{jN}$ is feasible to ($Q$) with probability ($1-\delta$). (5) If every ${\bf x}_{jN}$ ($j=1,\ldots, M$) is feasible, then, STOP; otherwise, increase $N$ and repeat steps (2) to (5). 

The pseudo codes of the lower and upper bound procedures are presented in Appendix A.

\subsection{Distributed and Parallel Implementation of the SAA Algorithm}
In our computational analysis, we solve $M$ replications of ($\hat{Q}$) to generate bounds for the true problem ($Q$). Solving ($\hat{Q}$) via the SAA Algorithm is time consuming because the MILP formulation ($\hat{Q}_{\alpha}$) is solved $\big\lceil\frac{\ln_2(\alpha_u - \alpha_l)}{\ln_2(\phi)}\big\rceil$ times; and each time with a different value of $\alpha$. Thus, we solve $M$ problems using  distributed  and parallel computing capabilities of the Plametto cluster, a high-performance computing resource of Clemson University. 

Via distributed computing, the computational problem is separated into several individual tasks, where each of task is completed on one or more nodes. Unlike the distributed computing, in parallel computing tasks are further divided into smaller tasks which are then solved  using multi-processors on a single computer \cite{tanenbaum2007distributed}. The literature presents several different parallel computing approaches, including bit-level, instruction-level, data, and task parallelism \cite{padua2011encyclopedia}. We use task parallelism, which execute different tasks across multiple threads at the same time.  

In the  distributed computing environment of the Palemtto cluster, let $N^c$ represent the computing nodes (computers) that communicate with one another to solve the $M$ different problems (replications), each assigned to a node. These problems are solved in parallel, that is, each problem is solved simultaneously and without any intersection or disturbance. Let  $P^t$ represent the total number of processors per node. In Appendix A we present the pseudo code of the distributed and parallel computing procedure we develop. 

\section{Case Study and Experimental Setup} \label{experiment}
\subsection{Case Study} \label{sec:CaseStudy}
We develop a case study to validate the models developed and derive managerial insights. The case study considers that a blend of corn stover 2-pass (2-P),  corn stover 3-pass (3-P), switchgrass and miscanthus is processed in a facility which uses biochemical conversion to produce renewable fuel. Figure \ref{flow chart} presents the flowchart of processes in this facility. This process is expected to  deliver a consistent feedstock which meets biochemical conversion specifications as designed by the National Renewable Energy Laboratory  with a total structural carbohydrate of 59.1\% and a moisture content of 20\% \cite{davis2013process}.  We utilize historical data collected at INL's PDU during the recent few years. Table \ref{tab:inputdata} in Appendix summarizes technical, process related data for switchgrass. 
 
\begin{figure*}[htb]
    \centering
    \includegraphics[width=0.7\columnwidth]{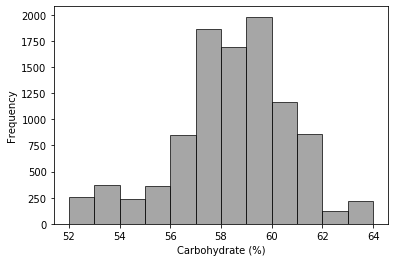}
    \includegraphics[width=0.7\columnwidth]{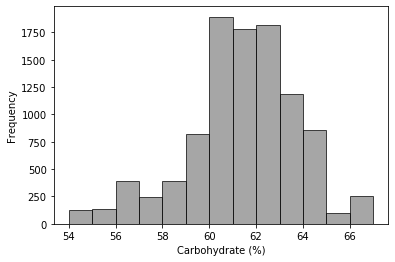}
    \caption{Empirical distribution of corn stover 3-P and corn stover 2-P.}
    \label{C3C2}
\end{figure*}

\begin{figure*}[htb]
    \centering
    \includegraphics[width=0.7\columnwidth]{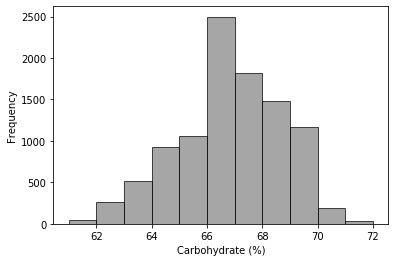}
    \includegraphics[width=0.7\columnwidth]{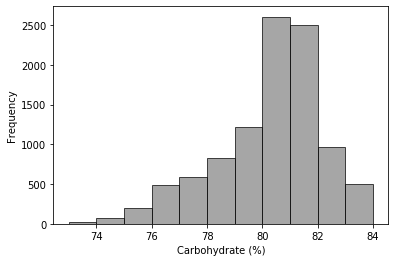}
    \caption{Empirical distribution of swichgrass and miscanthus.}
    \label{SM}
\end{figure*}

Figures \ref{C3C2} and \ref{SM} represent the empirical distribution of carbohydrate content for each feedstock. In order to develop these empirical distributions we used historical collected at INL, data collected at DuPont \cite{JR21}, and data reported by \cite{Owens2016RegionalFP}.     



\begin{table}[htbp]
  \centering
  \caption{Input Data: Biomass Blend}
    \scalebox{1}{\begin{tabular}{l|ccc|c|c} \toprule
    & \multicolumn{4}{c|}{\bf Nr. of Bales}   &\\
    \cline{2-5}
    \multicolumn{1}{c|}{\bf Feedstock}	& \multicolumn{3}{c|}{\bf Moisture Level}    & {\bf Total}  & {\bf Avg. Carb.}\\
    \multicolumn{1}{c|}{\bf Type}	& {\bf	L}    & {\bf M}  & {\bf H}    &   & {\bf Cont. (\%)}  \\   \midrule
Corn Stover (3-P)	&	6	&	10	&   4    &   20    &   57.4     \\
Corn stover (2-P)	&	12	&	20	&   8    &    40   &   60.3    \\
Switchgrass	&	3	&	5	&   2    &   10    &   66.6      \\
Miscanthus	&	3	&	5	&   2    &   10    &   81.7     \\
\bottomrule
\end{tabular}}%
  \label{Number of BalesS}%
\end{table}%

Table \ref{Number of BalesS} summarizes the data we use in the case study. We consider that 80 bales of biomass are to be processed in this system. For each feedstock type we report the total number of bales, the number of bales per moisture level, and the average carbohydrate content per bale.  Table \ref{Sequence setting} presents the sequences of bales we use. The performance of the system is impacted by the sequence, thus, we solve 6 problems, each using a different sequence. In these sequences we use letters H, L, M to represent high, low and medium moisture content correspondingly. We use letters C2, C3, S and M to represent corn stover 2-P, stover 3-P, switchgrass and miscanthus correspondingly. The numbers represent the number of bales. We also use a  random sequence, and a short sequence. The short sequence is presented in the Appendix.   

\begin{table}[htbp]
  \centering
  \caption{Problem Definitions}
    \scalebox{1}{ \begin{tabular}{c|l|l} \toprule
     {\bf Problem }	& \multicolumn{2}{c}{\bf Sequence Pattern}   \\  \cline{2-3}
    {\bf  Nr.}	& {\bf Feedstock-Based}  & {\bf Moisture-Based}   \\ \midrule
1	&	10S-40C2-10M-20C3	& 3L-3M-2H\\
2	&	10S-20C2-10M-20C2-20C3	& 3L-3M-2H\\
3	&	Short sequence	& 3L-3M-2H\\
4	&	Short sequence	& 5M-3L-2H\\
5   &   Short sequence   & 2H-3L-5M  \\
6   &   Random sequence & 3L-3M-2H\\
\bottomrule
\end{tabular}}%
  \label{Sequence setting}%
\end{table}%
 
\subsection{Experimental Setup} 
The algorithms proposed are coded in Python \cite{van1995python}. Gurobi Optimizer is used to solve the MILPs \cite{grb2010}. The parallel implementation of SAA Algorithm was conducted in the Plametto cluster of Clemson University. 

Tables \ref{Results0.05E} and \ref{Results0.1E}  summarize the results from our computational effort to identify a sample size $N$ and $\hat{\gamma}$ for which the solution of model ($\hat{Q}$) is feasible for model ($Q$), with high confidence level. We follow the procedure described in Section \ref{Sec:UpperLowerBd}. In these experiments, $M = 10,$ and $\gamma \geq \hat{\gamma}.$ We only report the results of solving problem 3 for different values of $N$ and $\hat{\gamma}$. The confidence level is $(1- \delta)*100 = 99\%.$ 

Based on our numerical results, when $\gamma > \hat{\gamma}$ (Table \ref{Results0.05E}) and $N$ = 400 or $N$ = 600, each of $M$ replications gave a feasible solution. The corresponding maximum solution risk is smaller then $\gamma.$ When $\gamma = \hat{\gamma}$ (Table \ref{Results0.1E}) and $N$ = 2,500, each of $M$ replication gave a feasible solution. The corresponding maximum solution risk is smaller then $\gamma.$ Notice that the computational efforts  to find feasible solutions when $\gamma = \hat{\gamma}$ are higher then in the case when $\gamma > \hat{\gamma}$. The quality of the optimal solutions for $\gamma = \hat{\gamma}$ is slightly better.

Table \ref{ResultsLB} summarizes the lower bounds generated for problem 3 using the procedure outlined in  Section \ref{Sec:UpperLowerBd}. Based on these results, the upper bounds generated when $\gamma = \hat{\gamma}$ and $N$ = 2,500 (Table \ref{Results0.1E}) are optimal at 99\% confidence level. The upper bounds generated when $\gamma > \hat{\gamma}$ and $N$ = 400 (Table \ref{Results0.05E}) are 0.3 to 0.5\% away from optimality at 99\% confidence level. Due to the high computational efforts to find feasible solutions when $\gamma = \hat{\gamma}$ (with only marginal improvement to solution quality), in the rest of our experiments we use $N = 400$, $M = 10$, and $\hat{\gamma}$ = 0.05. 
\begin{table}[htbp]
  \centering
  \caption{Upper Bounds: Problem 3,  \\$\gamma=0.10$, $\hat{\gamma}=0.05$}
    \scalebox{0.8}{ \begin{tabular}{c|ccc|ccc|c|c} \toprule
    {\bf Sample}	& \multicolumn{3}{c|}{\bf Process Time(hrs)}  & \multicolumn{3}{c|}{\bf Solution Risk} &  {\bf No. of}  &   \multicolumn{1}{c}{\bf Avg. Run}\\
{\bf Size ($N$)}	&   {\bf Min}    & {\bf Max}    &   {\bf Avg}   & {\bf Min}    & {\bf Max}    &   {\bf Avg}    &   {\bf Feasible}	&  {\bf Time (s)}  \\ \midrule
50	&	8.97	&	9.00	&   8.99    &      0.000   &  0.132   &   0.021    &   7    &      739    \\
100	&	8.95	&	8.98	&   8.97    &      0.000    &   0.141    &   0.016    &   9    &      1,010    \\
150	&	8.95	&	8.98	&   8.97    &     0.000    &   0.123    &   0.016    &   6    &      1,059    \\
200	&	8.93	&	8.98	&   8.95    &    0.000    &   0.118    &   0.016    &   7    &    1,208    \\
250	&	8.93	&	8.98	&   8.96    &     0.000    &   0.120    &   0.014    &   8    &     1,241    \\
300	&	8.97	&	8.98	&   8.97    &    0.000    &   0.116    &   0.012    &   9    &     1,353    \\
400	&	8.93	&	8.98	&   8.96    &    0.000    &   0.097    &   0.013    &   10    &     1,646    \\
600	&	8.93	&	8.97	&   8.96    &    0.000    &   0.097    &   0.013    &   10    &       2,176   \\
\bottomrule
\end{tabular}}%
\label{Results0.05E}%
\end{table}%

\begin{table}[htbp]
  \centering
  \caption{Upper Bounds: Problem 3,  $\gamma=\hat{\gamma}=0.1$}
    \scalebox{0.8}{ \begin{tabular}{c|ccc|ccc|c|c} \toprule
    {\bf Sample}	& \multicolumn{3}{c|}{\bf Process Time(hrs)}  & \multicolumn{3}{c|}{\bf Solution Risk} &  {\bf No. of}  &   \multicolumn{1}{c}{\bf Avg. Run}\\
{\bf Size ($N$)}	&   {\bf Min}    & {\bf Max}    &   {\bf Avg}   & {\bf Min}    & {\bf Max}    &   {\bf Avg}    &   {\bf Feasible}	&  {\bf Time (s)}  \\ \midrule
1000	&	8.93	&	8.93	&   8.93    &    0.000   &  0.105   &   0.013    &   7    &     2,744    \\
2500	&	8.93	&	8.93	&   8.93    &  0.000    &   0.099    &   0.013    &   10    &       6,390    \\
5000	&	8.93	&	8.93	&   8.93    &    0.000    &   0.100    &   0.013    &   9    &     12,715    \\
\bottomrule
\end{tabular}}%
\label{Results0.1E}%
\end{table}%

\begin{table}[htbp]
  \centering
  \caption{Lower Bounds: Problem 3, $\gamma=\hat{\gamma}=0.1$}
    \scalebox{0.8}{ \begin{tabular}{c|c|ccc|c} \toprule
    {\bf Sample}	& {\bf Lower}  & \multicolumn{3}{c}{\bf Gap(\%)}  & {\bf Ave. Running}\\
{\bf Size ($N$)}	&   {\bf Bound(hrs)}    &   {\bf Min}   &   {\bf Max}  & {\bf Ave}    &   {\bf Time (sec)}    \\ \midrule
1000	&	8.93	&	0.00	&   0.00    &   0.00    &   2,744   \\
2500	&	8.93	&	0.00	&   0.00    &   0.00    &   6,390    \\
5000	&	8.93	&	0.00	&   0.00    &   0.00    &   12,715   \\
\bottomrule
\end{tabular}}%
\label{ResultsLB}%
\end{table}%
 
\section{Analysis of Results and Observations} \label{result}
Tables \ref{Result1595} to \ref{Results3090} summarize experimental results when solving problems 1 through 6 for different values of $q$, $\underline{q}$ and $\tau.$ Recall that $\underline{q}$*100 and $q$*100 represent the minimum and maximum utilization of the reactor. The value of $\tau$ represents the frequency by which the chance constraint is satisfied. While constraint \eqref{eq:blendingf} enforces that the chance constraint is met within the time horizon $\mathcal{T};$  constraints \eqref{eq:blendingf_tau} enforce that the chance constraint is met every $\tau$ periods. In our experiments $\tau$ takes values 15 and 30 minutes. 
\begin{align}
Pr \left( \sum_{i \in \mathcal{I}_{r}} \sum_{s \in \mathcal{S}} \sum_{b \in \mathcal{B}}  (f_{b} - f^*) X_{isbt} \geq 0 \right ) & \geq 1 - \gamma \nonumber\\
\forall t = \tau, 2\tau, \ldots, \mathcal{T}. \label{eq:blendingf_tau} 
\end{align}

In our numerical analysis we report the value of the \emph{coefficient of variation,} which measures the variability of biomass flow to the reactor, and is calculated as the standard deviation over $t\in\mathcal{T}$ of: \[\frac{\sum_{i\in\mathcal{I}_r}\sum_{s\in\mathcal{S}}\sum_{b\in\mathcal{B}}X_{isbt}}{\sum_{i\in\mathcal{I}_r}\sum_{s\in\mathcal{S}}\sum_{b\in\mathcal{B}}\sum_{t\in\mathcal{T}}X_{isbt}}.\]

\noindent {\bf Observation 1:} \emph{The processing rate of the reactor decreases,  and the processing time, the maximum inventory, and the coefficient of variation increase as the requirements for meeting carbohydrate content become rigid.} The experiments summarized in Table \ref{Result1595} considers $\tau = 15$ and Table \ref{Results3095} considers $\tau = 30.$ When $\tau = 30$, the average processing time of problems 1 to 5 is 2.2\%  shorter, reactor's processing rate is 2.5\% higher, maximum inventory level is 6.5\% lower, coefficient of variation is 14.9\% lower.  A similar observation is made when comparing the results of Tables \ref{Results1590} and \ref{Results3090}. When $\tau = 30$, the average processing time of problems 1 to 5 is 1.6\%  shorter, reactor's processing rate is 1.6\% higher, maximum inventory level is 1.4\% lower, coefficient of variation is 4.5\% lower. Furthermore, problem 6, which represents the random sequence, could not find solutions which satisfied the chance constraints for $\tau = 15$ and $\tau = 30$ minutes.

\begin{table}[htbp]
  \centering
  \caption{Summary of Results: $\tau=15$ min, $\underline{q}=0.90, q = 0.95$}
    \scalebox{0.65}{ \begin{tabular}{c|c|c|ccc|c} \toprule
    {\bf Problem}	& {\bf Process Time}  & {\bf Reactor Rate}  & {\bf Flow to} &   {\bf Ave. Inv.}    &   {\bf Max. Inv.}    & {\bf Coefficient}   \\
{\bf Nr.}	&   {\bf (hrs)}    &   {\bf (mg/h)}   &   {\bf Reactor (mg)}  & {\bf (mg)}    & {\bf (mg)}    &   {\bf Variation}   \\ \midrule
1 & 13.62 & 2.26  & 30.75 & 1.12  & 4.54  & 0.0294 \\
2 & 11.70 & 2.63  & 30.75 & 0.83  & 4.18  & 0.0454 \\
3 & 8.96  & 3.43  & 30.75 & 0.15  & 0.87  & 0.0525 \\
4 & 9.30  & 3.31  & 30.75 & 0.17  & 0.98  & 0.0223 \\
5 & 13.25 & 2.32  & 30.75 & 0.98  & 4.22  & 0.0184 \\
6 & Infeasible &       &       &       &       &       \\
\bottomrule
\end{tabular}}%
  \label{Result1595}%
\end{table}%
\begin{table}[htbp]
  \centering
  \caption{Summary of Results: $\tau=30$ min, $\underline{q}=0.90, q = 0.95$}
    \scalebox{0.65}{ \begin{tabular}{c|c|c|ccc|c} \toprule
    {\bf Problem}	& {\bf Process Time}  & {\bf Reactor Rate}  & {\bf Flow to} &   {\bf Ave. Inv.}    &   {\bf Max. Inv.}    & {\bf Coefficient}   \\
{\bf Nr.}	&   {\bf (hrs)}    &   {\bf (mg/h)}   &   {\bf Reactor (mg)}  & {\bf (mg)}    & {\bf (mg)}    &   {\bf Variation}   \\ \midrule
1 & 13.62 & 2.26  & 30.75 & 1.12  & 4.54  & 0.0285 \\
2 & 10.78 & 2.86  & 30.75 & 0.63  & 3.31  & 0.0460 \\
3 & 8.96  & 3.43  & 30.75 & 0.14  & 0.83  & 0.0524 \\
4 & 8.96  & 3.43  & 30.75 & 0.13  & 0.87  & 0.0524 \\
5 & 13.25 & 2.32  & 30.75 & 0.98  & 4.28  & 0.0181 \\
6 & 15.15 & 2.04  & 30.75 & 1.31  & 4.53  & 0.0239 \\
\bottomrule
\end{tabular}}%
  \label{Results3095}%
\end{table}%
\begin{table}[htbp]
  \centering
  \caption{Summary of Results: $\tau=15$ min, $\underline{q}=0.80, q = 0.90$}
    \scalebox{0.65}{ \begin{tabular}{c|c|c|ccc|c} \toprule
    {\bf Problem}	& {\bf Process Time}  & {\bf Reactor Rate}  & {\bf Flow to} &   {\bf Ave. Inv.}    &   {\bf Max. Inv.}    & {\bf Coefficient}   \\
{\bf Nr.}	&   {\bf (hrs)}    &   {\bf (mg/h)}   &   {\bf Reactor (mg)}  & {\bf (mg)}    & {\bf (mg)}    &   {\bf Variation}   \\ \midrule
1 & 13.45 & 2.29  & 30.75 & 1.12  & 4.54  & 0.0575 \\
2 & 10.70 & 2.88  & 30.75 & 0.68  & 3.41  & 0.0965 \\
3 & 8.93  & 3.44  & 30.75 & 0.19  & 1.23  & 0.1108 \\
4 & 8.96  & 3.43  & 30.75 & 0.17  & 1.08  & 0.0955 \\
5 & 11.82 & 2.60  & 30.75 & 0.74  & 3.20  & 0.0363 \\
6 & Infeasible &       &       &       &       &       \\
\bottomrule
\end{tabular}}%
  \label{Results1590}%
\end{table}%
\begin{table}[htbp]
  \centering
  \caption{Summary of Results: $\tau=30$ min, $\underline{q}=0.80, q = 0.90$}
    \scalebox{0.65}{ \begin{tabular}{c|c|c|ccc|c} \toprule
    {\bf Problem}	& {\bf Process Time}  & {\bf Reactor Rate}  & {\bf Flow to} &   {\bf Ave. Inv.}    &   {\bf Max. Inv.}    & {\bf Coefficient}   \\
{\bf Nr.}	&   {\bf (hrs)}    &   {\bf (mg/h)}   &   {\bf Reactor (mg)}  & {\bf (mg)}    & {\bf (mg)}    &   {\bf Variation}  \\ \midrule
1 & 13.42 & 2.29  & 30.75 & 1.11  & 4.54  & 0.0583 \\
2 & 9.90  & 3.11  & 30.75 & 0.49  & 2.78  & 0.0985 \\
3 & 8.93  & 3.44  & 30.75 & 0.19  & 1.22  & 0.1108 \\
4 & 8.93  & 3.44  & 30.75 & 0.17  & 1.23  & 0.1107 \\
5 & 11.82 & 2.60  & 30.75 & 0.74  & 3.49  & 0.0363 \\
6 & 14.00 & 2.22  & 30.75 & 1.17  & 4.45  & 0.0524 \\
\bottomrule
\end{tabular}}%
  \label{Results3090}%
\end{table}%

\noindent {\bf Observation 2:} \emph{Deterministic modeling of the system leads to underestimation of the processing rate of the reactor, processing time, inventory and coefficient of variation.} Tables \ref{ResultDeter1595} and \ref{ResultsDeter3095} summarize the results from solving the deterministic equivalent (the mean value problem) of problems 1 to 6. We obtain the deterministic equivalent model by substituting constraint \eqref{eq:blendingf} with constraints \eqref{eq:blendingf_deta}. In \eqref{eq:blendingf_deta},  $\bar{f}_{b}$ represents the average carbohydrate rate of feedstock $b$. Since \eqref{eq:blendingf_deta} are  deterministic constraints, we do not need to use the SAA method.
\begin{align}
\sum_{i \in \mathcal{I}_{r}} \sum_{s \in \mathcal{S}}\sum_{b \in \mathcal{B}}  (\bar{f}_{b} - f^*) X_{isbt} \geq 0  \hspace{0.3in} \forall t = \tau, 2\tau, \ldots, \mathcal{T}. \label{eq:blendingf_deta} 
\end{align}

\begin{table}[htbp]
  \centering
  \caption{Summary of Deterministic Results: $\tau=15$ min, $\underline{q}=0.90, q = 0.95$}
    \scalebox{0.65}{ \begin{tabular}{c|c|c|ccc|c} \toprule
    {\bf Problem}	& {\bf Process Time}  & {\bf Reactor Rate}  & {\bf Flow to} &   {\bf Avg Inv.}    &   {\bf Max Inv}    & {\bf Coefficient}   \\
{\bf Nr.}	&   {\bf (hrs)}    &   {\bf (mg/h)}   &   {\bf Reactor (mg)}  & {\bf (mg)}    & {\bf (mg)}    &   {\bf Variation}   \\ \midrule
1 & 9.00  & 3.42  & 30.75 & 0.14  & 1.09  & 0.0523 \\
2 & 9.00  & 3.42  & 30.75 & 0.15  & 0.98  & 0.0522 \\
3 & 8.93  & 3.44  & 30.75 & 0.10  & 0.72  & 0.0522 \\
4 & 8.93  & 3.44  & 30.75 & 0.11  & 0.92  & 0.0525 \\
5 & 13.25 & 2.32  & 30.75 & 0.98  & 4.54  & 0.0184 \\
6 & 13.23 & 2.32  & 30.75 & 1.04  & 3.31  & 0.0118 \\
\bottomrule
\end{tabular}}%
  \label{ResultDeter1595}%
\end{table}%

\begin{table}[htbp]
  \centering
  \caption{Summary of Deterministic Results: $\tau=30$ min, $\underline{q}=0.90, q = 0.95$}
    \scalebox{0.65}{ \begin{tabular}{c|c|c|ccc|c} \toprule
    {\bf Problem}	& {\bf Process Time}  & {\bf Reactor Rate}  & {\bf Flow to} &   {\bf Avg Inv}    &   {\bf Max Inv}    & {\bf Coefficient}   \\
{\bf Nr.}	&   {\bf (hrs)}    &   {\bf (mg/h)}   &   {\bf Reactor (mg)}  & {\bf (mg)}    & {\bf (mg)}    &   {\bf Variation}   \\ \midrule
1 & 8.98  & 3.42  & 30.75 & 0.15  & 1.25  & 0.0525 \\
2 & 8.98  & 3.42  & 30.75 & 0.15  & 1.05  & 0.0526 \\
3 & 8.93  & 3.44  & 30.75 & 0.10  & 0.69  & 0.0523 \\
4 & 8.93  & 3.44  & 30.75 & 0.11  & 0.91  & 0.0525 \\
5 & 13.25 & 2.32  & 30.75 & 0.98  & 2.96  & 0.0183 \\
6 & 13.23 & 2.32  & 30.75 & 1.05  & 4.52  & 0.0116 \\
\bottomrule
\end{tabular}}%
  \label{ResultsDeter3095}%
\end{table}%
By comparing the results of Tables \ref{Result1595} and \ref{Results3095} with Tables \ref{ResultDeter1595} and \ref{ResultsDeter3095}, respectively, we observe that the processing time of the deterministic model is 9.67\% shorter, reactor's processing rate is 11.93\% higher, maximum inventory level is 89.28\% lower and the coefficient of variation is 78.11\% lower then the corresponding results of the stochastic model. 

\noindent {\bf Observation 3:} \emph{Short sequences that begin by processing bales with low and medium moisture level perform best. These sequences lead to robust solutions.} Based on the results in Tables \ref{Result1595} through \ref{Results3090}, problems 3 and 4 perform the best, and problem 3 outperforms problem 4. In both problems we use the short feedstock-based sequence. For problem 3 the moisture-based sequence begins by processing low moisture bales, followed by medium  and high moisture level. Such sequences allow the system to build up inventory which is used when processing high moisture bales. Problems 1 and 2 perform worst then 3, 4 and 5. In these problems, we process every bale of a particular feedstock at a time. Such a strategy leads to higher inventory level and longer processing times. In many of the problem instances we solved, the Random sequence could not find a feasible solution. In the cases when a feasible solution is found, its quality is the worst than problems 1 to 5. This points to the importance of sequencing bales to streamline biomass preprocessing.  

Comparing the results of problem 3 in Tables \ref{Result1595} through \ref{Results3090} with Tables \ref{ResultDeter1595} and \ref{ResultsDeter3095} one can see that the performance of problem 3 did not change. This indicates that problem 3 is robust to variations in problem parameters. 

In our numerical analysis we assume that the time unit $t=1$ minute ($t\in\mathcal{T}$). We resolve Problem 3 considering $t =0.25$ minutes. The results are summarized in Table \ref{Results0.25}. Using smaller values of $t$ increase the size of the problem and the running time of our proposed algorithms. However, doing this leads to an 8.3\% reduction of processing time of the system, 9\% increase of processing rate of the reactor, and 50\% lower maximum inventory.   
\begin{table}[H]
  \centering
  \caption{Summary of Stochastic Results: $\tau=15$, $\underline{q}=0.90, q = 0.95$, $t=0.25$min}
    \scalebox{0.65}{ \begin{tabular}{c|c|c|ccc|c} \toprule
    {\bf Problem}	& {\bf Process Time}  & {\bf Reactor Rate}  & {\bf Flow to} &   {\bf Ave. Inv.}    &   {\bf Max. Inv.}    & {\bf Coefficient}   \\
{\bf Nr.}	&   {\bf (hrs)}    &   {\bf (mg/h)}   &   {\bf Reactor (mg)}  & {\bf (mg)}    & {\bf (mg)}    &   {\bf Variation}   \\ \midrule
3 & 8.19  & 3.75  & 30.75 & 0.15  & 0.82  & 0.0524 \\
\bottomrule
\end{tabular}}%
  \label{Results0.25}%
\end{table}%
\noindent {\bf Observation 4:} \emph{Increasing reactor's reliability requires the use of high quality biomass and additional inventory.} 
Tables \ref{Result1595} and \ref{Results3095} assume $q = 0.95$ and $\underline{q} = 0.90$, and Tables \ref{Results1590} and \ref{Results3090} assume $q = 0.90$ and $\underline{q} = 0.80$. The 5.56\% increase of the average reactor utilization ($q$) from 0.90 to 0.95 led to  4.92\% longer processing time, 2.44\% decrease of reactor's processing rate, 6.64\% increase of maximum inventory. We conclude that  a system can achieve a higher utilization rate of the reactor via additional inventory level  and high biomass quality.   

\noindent {\bf Observation 5:} \emph{Chance constraint optimization modeling of the system identifies robust and realistic solutions to our problem.} We resolve model ($Q$) as a stochastic optimization model. To accomplish this, we set $\hat{\gamma} = 0.00$, thus,  constraint \eqref{eq:blendingf} becomes. 
\begin{align}
Pr \left( \sum_{i \in \mathcal{I}_{r}} \sum_{s \in \mathcal{S}} \sum_{b \in \mathcal{B}}  (f_{b} - f^*) X_{isbt} \geq 0 \right ) & \geq 1 \nonumber\\ \forall t = \tau, 2\tau, \ldots, \mathcal{T}. \label{eq:blendingf_det} 
\end{align}
Since in \eqref{eq:blendingf_det} $P() \geq 1$,  these constraints are to be satisfied all the time. Thus, constraints \eqref{eq:blendingf_det} and \eqref{eq:blendingf_det2} are equivalent. Thus, the solution of ($\hat{Q}$) for $\hat{\gamma} = 0.00$ is the corresponding stochastic optimization solution of ($\hat{Q}$).
\begin{align}
\sum_{i \in \mathcal{I}_{r}} \sum_{s \in \mathcal{S}} \sum_{b \in \mathcal{B}}  (f_{b} - f^*) X_{isbt} \geq 0  \nonumber\\ \forall t = \tau, 2\tau, \ldots, \mathcal{T}. \label{eq:blendingf_det2} 
\end{align}

\begin{table}[htbp]
  \centering
  \caption{Summary of Stochastic Results: $\tau=15$ min, $\underline{q}=0.90, q = 0.95$}
    \scalebox{0.6}{ \begin{tabular}{c|c|c|ccc|c|c} \toprule
    {\bf Problem}	& {\bf Process Time}  & {\bf Reactor Rate}  & {\bf Flow to} &   {\bf Avg Inv}    &   {\bf Max Inv}    & {\bf Coefficient}   &   {\bf Nr. of}   \\
{\bf Nr.}	&   {\bf (hrs)}    &   {\bf (mg/h)}   &   {\bf Reactor (mg)}  & {\bf (mg)}    & {\bf (mg)}    &   {\bf Variation}   &   {\bf Feasible} \\ \midrule
1 & Infeasible &       &       &       &       &       & 0 \\
2 & Infeasible &       &       &       &       &       & 0 \\
3 & 10.21 & 3.05  & 30.75 & 0.47  & 2.08  & 0.0475 & 7 \\
4 & 12.64 & 2.47  & 30.75 & 0.92  & 4.12  & 0.0276 & 7 \\
5 & 13.41 & 2.30  & 30.75 & 1.01  & 4.30  & 0.0191 & 7 \\
6 & Infeasible &       &       &       &       &       & 0 \\
\bottomrule
\end{tabular}}%
  \label{ResultWithout1595}%
\end{table}%

\begin{table}[htbp]
  \centering
  \caption{Summary of Stochastic Results: $\tau=30$ min, $\underline{q}=0.90, q = 0.95$}
    \scalebox{0.6}{ \begin{tabular}{c|c|c|ccc|c|c} \toprule
    {\bf Problem}	& {\bf Process Time}  & {\bf Reactor Rate}  & {\bf Flow to} &   {\bf Avg Inv}    &   {\bf Max Inv}    & {\bf Coefficient}   &   {\bf Nr. of} \\
{\bf Nr.}	&   {\bf (hrs)}    &   {\bf (mg/h)}   &   {\bf Reactor (mg)}  & {\bf (mg)}    & {\bf (mg)}    &   {\bf Variation}   &   {\bf Feasible}  \\ \midrule
1 & Infeasible &       &       &       &       &      & 0 \\
2 & 14.17 & 2.17  & 30.75 & 1.25  & 4.54  & 0.0517 & 1 \\
3 & 10.18 & 3.06  & 30.75 & 0.46  & 1.99  & 0.0477 & 7 \\
4 & 10.43 & 2.99  & 30.75 & 0.48  & 2.15  & 0.0392 & 7 \\
5 & 13.40 & 2.30  & 30.75 & 1.01  & 4.28  & 0.0199 7 \\
6 & Infeasible &       &       &       &       &       & 0 \\
\bottomrule
\end{tabular}}%
  \label{ResultsWithout3095}%
\end{table}%

Tables \ref{ResultWithout1595} and \ref{ResultsWithout3095} summarize the results of the stochastic programming formulation of ($\hat{Q}$). When $\tau = 15,$ we did not find a feasible solution in 3 out of 6 problems solved. Of the problems for which we found feasible solutions, 7 out of the 10 replications were feasible. A similar observation is made when $\tau = 30.$ Furthermore, comparing the results of Tables \ref{Result1595} and \ref{Results3095} with Tables \ref{ResultWithout1595} and \ref{ResultsWithout3095}, we observe that the stochastic model leads to 19.21\% increase in processing time, 7.26\% decrease in reactor's processing rate, 17.10\% increase of maximum inventory, and 1.20\% decrease in the coefficient of variation.  

\section{Managerial Insights} \label{insight}
\subsection{Feasibility of Large-Scale Implementation of the Proposed Models}
Successful large-scale implementation of the models developed depends on a number of factors, such as, ($i$) standardizing the format of bales to facilitate transportation and storage; ($ii$) use of Radio Frequency Identification (RFID) technology to record and keep track of bale properties from harvest to the biorefinery; ($iii$) use of sensing and real time monitoring of material attributes, such as, moisture content since it changes with time; ($iv$) use of automated material handling equipment to efficiently move large quantities of biomass; ($v$) use of automated process control to implement the sequencing algorithms we propose.   

The focus of this research is on identifying biomass blends to meet process requirements within a biorefinery, given the inventory of biomass bales. The scope of this model could be extended to optimize the performance of the overall supply chain. Such a model would focus on identifying how many bales of different biomass feedstock to purchase given biomass availability in the region, biomass price and quality. The results of this model can be used to inform the design of long-term contracts among farmers and the biorefinery.

\subsection{Simple Rules to Sequence Bales}
The results of our numerical analysis point to the importance of sequencing bales based on moisture level and carbohydrate content. Although in this paper we present the results from solving 6 problems, each representing a different sequence, we tested numerically a number of different sequences. We conducted numerical tests, rather then solving an optimization model, for two main reasons. First, integrating a sequencing model within our chance constraint optimization would increase the computational challenge of solving the problem. Second, since our model considers only 3 moisture levels and 4 feedstocks, we were able to identify numerically the characteristics of sequences which work best. Below we provide  simple rules that practitioners can use to sort bales based on moisture level and carbohydrate contents.       

\noindent {\bf Rule 1:} \emph{Sort Bales Based on Moisture Level}: Let $\mathcal{N}_L, \mathcal{N}_M, \mathcal{N}_H$ represent the number of bales of low, medium and high moisture level. Let $\psi(\mathcal{N}_L, \mathcal{N}_M, \mathcal{N}_H)$ represent the greatest common divisor (GCD) of $\mathcal{N}_L, \mathcal{N}_M, \mathcal{N}_H$.  

\emph{The rule is}: Begin the sequence by processing $\frac{\mathcal{N}_L}{\psi()}$ low moisture bales, followed by $\frac{\mathcal{N}_M}{\psi()}$ medium moisture bales, followed by $\frac{\mathcal{N}_H}{\psi()}$ high moisture bales. Repeat this sequence $\psi()$ times to process every bale in the inventory.  

In the case when the GCD of $\mathcal{N}_L, \mathcal{N}_M, \mathcal{N}_H$ does not exist, one can adjust these number by decreasing them by 1 or 2 bales so that the resulting  $\mathcal{N}^{\prime}_L, \mathcal{N}_M^{\prime}, \mathcal{N}_H^{\prime}$ have a GCD. Next, the same Rule 1 applies for $\mathcal{N}^{\prime}_L + \mathcal{N}_M^{\prime} + \mathcal{N}_H^{\prime}$ bales. The remaining bales are processed at the very end. 

\underline{Example:} Let $\mathcal{N}_L = 24, \mathcal{N}_M = 40, \mathcal{N}_H = 16$ bales. Thus, $\psi(24, 40, 16) = 8.$ The sequence is $3L-5M-2H$, which repeats 8 times. 
 
\noindent {\bf Rule 2:} \emph{Sort Bales Based on Biomass Characteristics}: 
Let $\mathcal{N}_q$ represent the number of bales that meet carbohydrate requirements, and $\mathcal{N}_{nq}$ represent the number of bales that do not. 
 Let $\psi(\mathcal{N}_q,\mathcal{N}_{nq})$ be the GCD of $\mathcal{N}_q$ and $\mathcal{N}_{nq}$. 

\emph{The rule is}: Begin the sequence by processing $\frac{\mathcal{N}_q}{\psi}$ bales, followed by $\frac{\mathcal{N}_{nq}}{\psi}$ bales.  Repeat this sequence $\psi$ times to process every bale in the inventory.  

In the case when the GCD of $\mathcal{N}_q, \mathcal{N}_{nq}$  does not exist, one can reduce $\mathcal{N}_{nq}$ to $\mathcal{N}^{\prime}_{nq}$ so that the resulting  $\mathcal{N}_{q}, \mathcal{N}_{nq}^{\prime}$ have a GCD. Next, the same Rule 2 applies for $\mathcal{N}_{q}+ \mathcal{N}^{\prime}_{nq}$ bales. The remaining ($\mathcal{N}_{nq} - \mathcal{N}_{nq}^{\prime}$) bales are processed at the very end.

\underline{Example:} Let $\mathcal{N}_q = 20, \mathcal{N}_{nq} = 60$ bales. Thus, $\psi(20, 60) = 20.$ The sequence is $1q-3{nq}$, which repeats 20 times.

\noindent {\bf Rule 3:} Let $f_b^{\gamma}$ represent the value of carbohydrate content which corresponds to the $\gamma$-th percentile of biomass bales of feedstock $b$. Let ${S}_q = \{b\in\mathcal{B}|f_b^{\gamma} \geq f^* \}$ and ${S}_{nq} = \{b\in\mathcal{B}|f_b^{\gamma} < f^*\}$.  Let $\underline{f} = \min_{b\in S_q} f_b^{\gamma}$, and let $N(\underline{f})$ be the corresponding number of bales. Let $\overline{f} = \max_{b\in S_{nq}} f_b^{\gamma}$, and let $N(\overline{f})$ be the corresponding number of bales.  The rule is outlined in Algorithm \ref{Rule_3}.

\begin{algorithm}
	\caption{Rule 3 Procedure}
	\textbf{Initialize:} $S_q$,  $S_{nq}$,  $f^*$
	\begin{algorithmic}[1]
	    \State \textbf{Calculate}: $\underline{f}$, $\overline{f}$, $N(\underline{f})$, $N(\overline{f})$
	   \State \textbf{Calculate}: $\overline{d} = \underline{f} - f^*$ and $\underline{d} = f^* - \overline{f}$;
		\State \textbf{While} $\left(N(\underline{f}) > 0 \And N(\overline{f}) > 0\right)$ {\bf Do}
		    \State \hspace{0.25in} {\bf If } ($\overline{d} \geq \underline{d}$): 
		     \State \hspace{0.5in}  Sequence 2 bales of $S_q$ followed by $\lfloor\frac{\overline{d}}{\underline{d}}\rfloor$ bales of $S_{nq}$
		     \State \hspace{0.5in} Let $N(\underline{f}) = N(\underline{f}) - 2$, $N(\overline{f}) = N(\overline{f}) - \lfloor\frac{\overline{d}}{\underline{d}}\rfloor$
		    \State \hspace{0.25in} {\bf Else:}  Sequence $\lceil\frac{\underline{d}}{\overline{d}}\rceil + 1$ bales of $S_q$ followed by 1 bale of $S_{nq}$
		     \State \hspace{0.5in} Let $N(\underline{f}) = N(\underline{f}) - \lceil\frac{\underline{d}}{\overline{d}}\rceil - 1$, $N(\overline{f}) = N(\overline{f}) - 1$
    	\State \textbf{EndWhile}
    	\State \textbf{If} $N(\underline{f}) = 0$ 
    	\State \hspace{0.25in}  $S_q = S_q \setminus b$
    	\State \textbf{Else}  
    	\State \hspace{0.25in}  $S_{nq} = S_{nq} \setminus b$
    	\State \textbf{If} $\left(S_q \neq \emptyset \And S_{nq} \neq \emptyset\right)$ , {\bf GoTo} line 1, otherwise {\bf GoTo} line 15
		\State Add the remaining bales to the sequence; and return the sequence created.
	\end{algorithmic} \label{Rule_3}
\end{algorithm}

\noindent {\bf Rule 4:} \emph{Sort Bales Based on Moisture Level and Biomass Characteristics}: This rule builds upon rules 1, 2 and 3.  

\emph{The rule is:} ($i$) Sequence bales based on moisture level using  rule 1.  
($ii$) {Sequence bales based on biomass characteristics} using rules 2 or 3. ($iii$) Sequence bales in such a way that both sequences identified are followed simultaneously. 

\underline{Example:} Let $\mathcal{N}_L = 24, \mathcal{N}_M = 40, \mathcal{N}_H = 16$ bales and $\mathcal{N}_q = 20, \mathcal{N}_{nq} = 60$ bales. Following rule 1 we get $3L-5M-2H$. Following rule 2 we get $1q-3nq$. Let $S_q =\{\text{Miscanthus, Switchgrass}\}$ and  $S_{nq} =\{\text{Corn stover 2-P, Corn stover 3-P.}\}$  
\begin{figure}[htb]
    \centering
    \includegraphics[width=1\columnwidth]{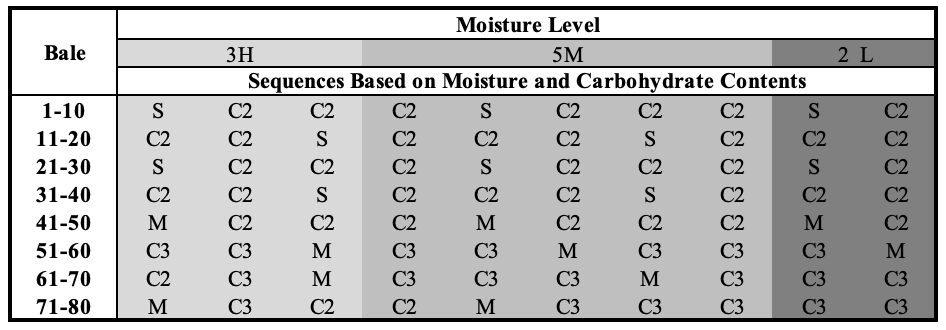}
    \caption{Sequence of problem 7.}
    \label{fig:sequence}
\end{figure}
Figure \ref{fig:sequence} demonstrates the sequence we created using the data in Section \ref{sec:CaseStudy}. This sequence is created using rule 4. We refer to this as problem 7. Notice that, both sequences are followed. The shade of the table changes with the moisture level to follow the pattern 3H-5M-2L. The values in the table change to represent the feedstock used, and these feedstocks  change based on 1q-3nq. Table \ref{ResultShort1595} summarizes the results of solving problem 7.  Comparing the  performance of problem 3 in Table \ref{Result1595} with the results of Table \ref{ResultShort1595}, it is clear that the performance of the sequence we develop using rule 4 is similar to the performance of the short sequence (problem 3).  
\begin{table}[htbp]
  \centering
  \caption{Rule 4: Summary of Results:
  $\tau=15$ min, $\underline{q}=0.90, q = 0.95$}
    \scalebox{0.65}{ \begin{tabular}{c|c|c|ccc|c} \toprule
    {\bf Problem}	& {\bf Process Time}  & {\bf Reactor Rate}  & {\bf Flow to} &   {\bf Ave. Inv.}    &   {\bf Max. Inv.}    & {\bf Coef. of}   \\
{\bf Nr.}	&   {\bf (hrs)}    &   {\bf (mg/h)}   &   {\bf Reactor (mg)}  & {\bf (mg)}    & {\bf (mg)}    &   {\bf Var.}   \\ \midrule
7 & 8.93  & 3.44  & 30.75 & 0.13  & 1.11  & 0.0514 \\
\bottomrule
\end{tabular}}%
  \label{ResultShort1595}%
\end{table}%

\section{Summary and Conclusions}
The objective of the proposed research is to develop process control strategies  to improve
the performance of a biorefinery, and ensure a continuous flow of biomass to the
reactor, while meeting the requirements of biochemical conversion process. The conversion process requires that the carbohydrate content of the biomass fed to the reactor is at least 59.1\%.  Satisfying this threshold is challenging because carbohydrate content is stochastic. We propose a stochastic optimization model to identify the infeed rate of biomass to the system, processing
speed of equipment and inventory level to ensure a continuous feeding of the reactor. We propose a sample average approximation (SAA) method to approximate the model, and develop an algorithm which solved the problem efficiently. The model is tested and validated using data from the Process Development Unit of Idaho National Laboratory.

We conduct extensive numerical analysis. The results of this analysis lead to several observations, such as, the robustness of the short bale sequences we propose. The corresponding process control  leads to improvements on processing time, processing rate of the reactor, inventory and coefficient of variations (in feeding biomass to the reactor). Finally, the implementation of the models developed requires the use of standardized bale format, Radio Frequency Identification technology,  sensing and real time monitoring of material attributes, automated material handling equipment, and automated process control. The scope of the model proposed can be extended to include the whole supply chain. These models can be used to identifying how many bales of different biomass feedstock to purchase given biomass availability in the region, biomass price and quality. Thus, the results of this model can be used to inform the design of long-term contracts among farmers and the biorefinery.
 
\bibliographystyle{plain}
\bibliography{Biomass_bibliography}
\section*{Acknowledgment}
The funding of this work was provided by the U.S. Department of Energy, Office of Energy Efficiency and Renewable Energy, Bioenergy Technologies Office under award Number DE- EE0008255 and Department of Energy Idaho Operations Office Contract No. DE-AC07-05ID14517. The authors greatly thank INL staff in the DOE Biomass Feedstock National User Facility for providing technical data for this study. Specifically, we acknowledge Neal A. Yancey and Jaya S. Tumuluru, who are the individual INL contributors to this research who provided helpful comments, data, and other forms of support for this analysis. The views expressed herein are those of the author only, and do not necessarily represent the views of DOE or the U.S. Government.

\begin{IEEEbiography}[{\includegraphics[width=1in,height=1.25in,clip,keepaspectratio]{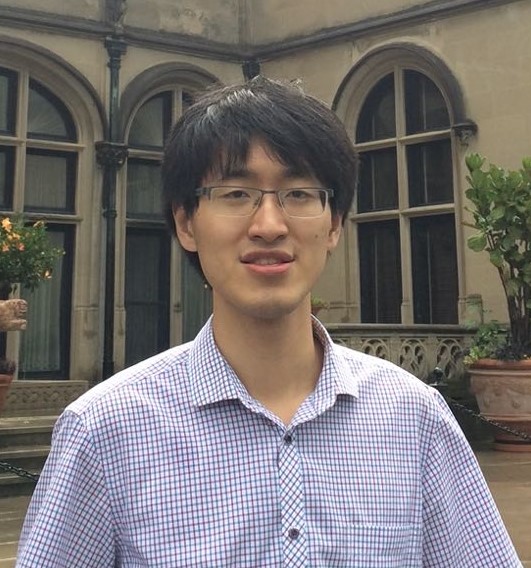}}]{Dr. Dahui (Kevin) Liu}
received his Ph.D in industrial engineering from Clemson University, Clemson, SC, USA, in 2020. Dr. Liu has worked for SAIC Motor Corporation Limited, Shanghai Aircraft Manufacturing Corporation Limited, and thyssenkrupp Uhde Chlorine Engineers (Shanghai) Corporation Limited. He is currently a Postdoctoral fellow in the Industrial Engineering Department, University of Arkansas. His research focuses in the area of operation research with applications in transportation and energy.
\end{IEEEbiography}
\vspace{-4in}
\begin{IEEEbiography}[{\includegraphics[width=1in,height=1.25in,clip,keepaspectratio]{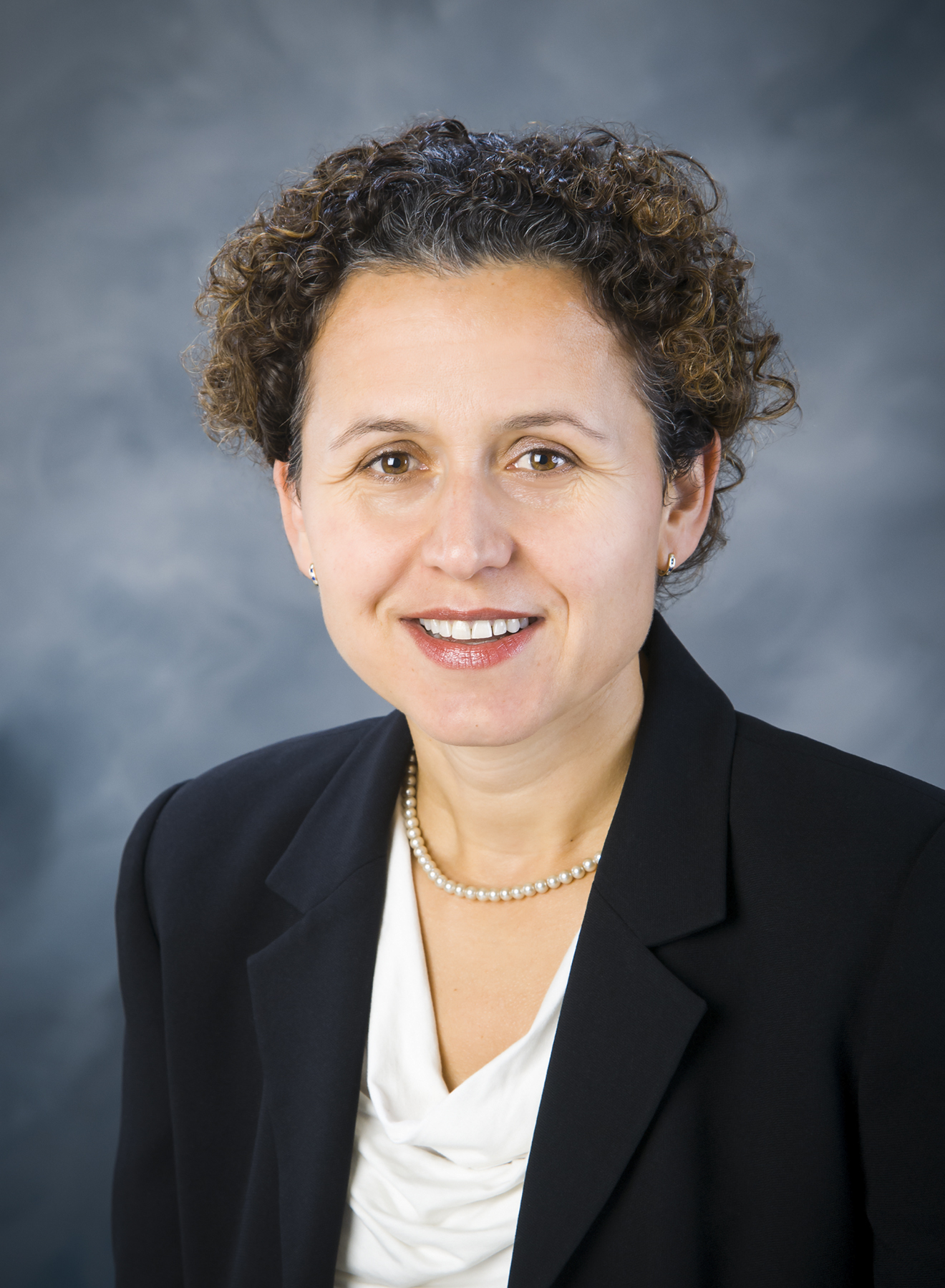}}] {Dr. Sandra D. Eksioglu} is the John M. and Marie G. Hefley Professor in Logistics and Entrepreneurship at the Department of Industrial Engineering, University of Arkansas. Dr. Eksioglu’s expertise is in the areas of operations research, network optimization, and algorithmic development with applications in transportation, logistics, and supply chain. 
\end{IEEEbiography}
\vspace{-4in}
\begin{IEEEbiography}[{\includegraphics[width=1in,height=1.25in,clip,keepaspectratio]{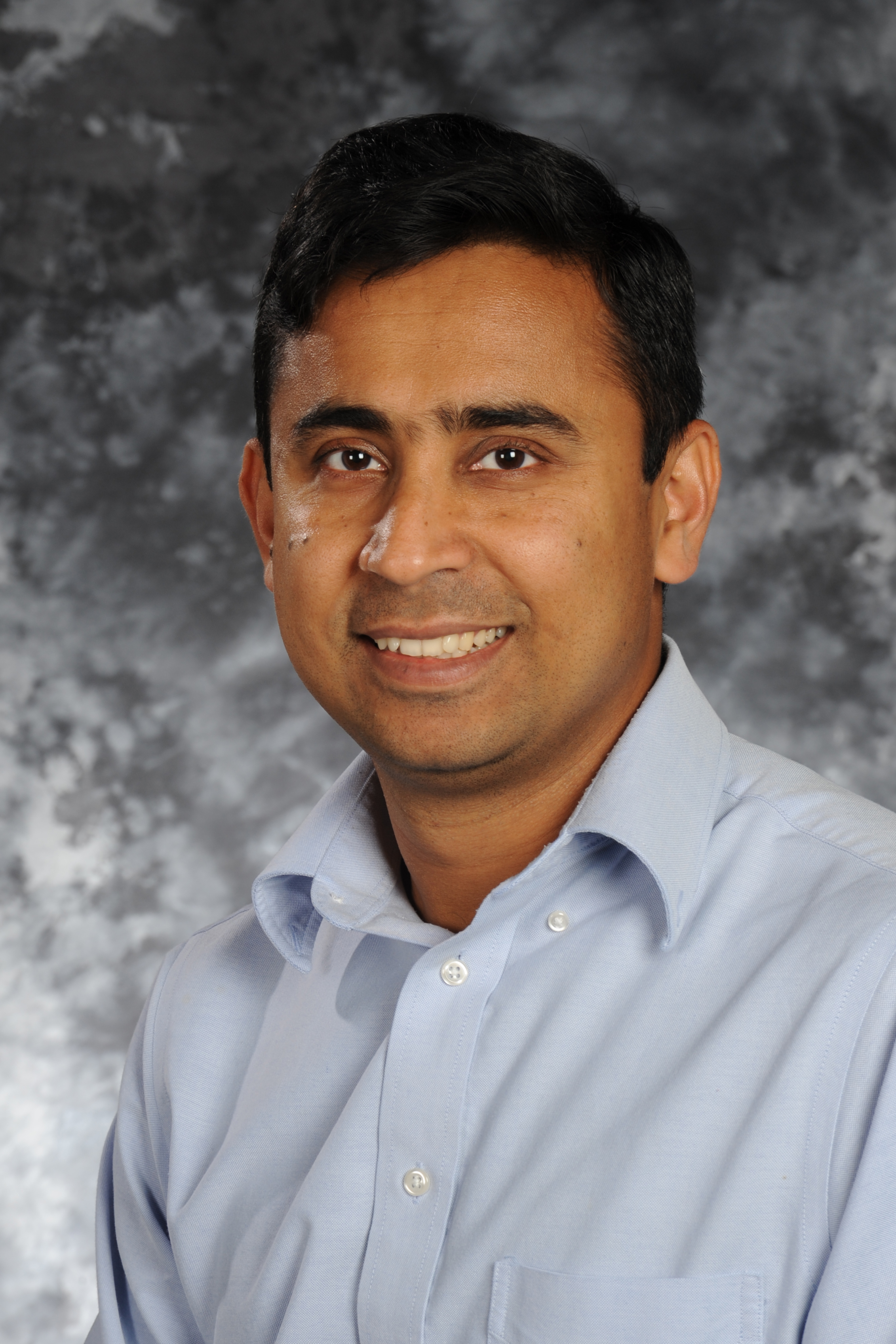}}] {Dr. Mohammad Roni} --``Roni" to his friends and colleagues -- is the Operations Research Lead of Idaho National Laboratory’s Software and Data Sciences department. In this position, he serves as principal investigator and technical lead on research projects sponsored by the U.S. Department of Energy. He is developing complex optimization models and algorithms applied to supply chain network analysis, resiliency optimization, transportation and infrastructure planning, and process control and design.
\end{IEEEbiography}

\clearpage
\appendices

\section{Algorithms}

\begin{algorithm}\label{Lower Bound Procedure}
	\caption{Lower Bound Procedure}
	\textbf{Initialize:} $\delta=0.01$; $\hat{\gamma} > \gamma$; $N\geq \frac{1}{2(\hat{\gamma} - \gamma)^2}\log{\big(\frac{1}{\delta}\big)}$; $M = 10$.
	\begin{algorithmic}[1]
	    \State $\underline{\nu} \leftarrow \infty$
		\For {$j = 1,\ldots,M$}
		    \State Sample $f_{bN}$ for $b\in\mathcal{B}$  
		    \State Solve ($\hat{Q}$) using SAA Algorithm to find $\hat{\nu}_{Nj}$
		    \If {$\hat{\nu}_{Nj} \leq \underline{\nu}$}
		        \State $\underline{\nu} \leftarrow \hat{\nu}_{Nj}$
		    \EndIf
		\EndFor
		\State Return $\underline{\nu}$
	\end{algorithmic} 
\end{algorithm} 

\begin{algorithm}\label{Upper Bound Procedure}
	\caption{Upper Bound Procedure}
	\textbf{Initialize:} $\delta = 0.01$; $C = 0$; $\hat{\gamma} < \gamma$; $N$; d; $N'=10000$; $M = 10$.
	\begin{algorithmic}[1]
	    \While{$C < M$}
	        \State $C \leftarrow 0$
		    \For {$j = 1,\ldots,M$}
		        \State Sample $f_{bN}$ for $b\in\mathcal{B}$  
		        \State Solve ($\hat{Q}$) using SAA Algorithm to find ${\bf x}_N$
		        \State Sample $f_{bN^{\prime}}$ for $b\in\mathcal{B}$  
		        \State Calculate $\bar{H}_{N^{\prime}}({\bf x}_{N},f)$ 
		        \State Calculate $\mathcal{U}_{N^{\prime}}({\bf x}_{N},f)$
		        \If {$\mathcal{U}_{N^{\prime}}({\bf x}_{N},f) \leq \gamma$}
		            \State ${\bf x}_N$ is a feasible to ($Q$)  
		            \State $C= C$ + 1
		        \Else
		            \State ${\bf x}_N$ is not feasible to ($Q$)
		        \EndIf
		    \EndFor
		    \State  $N = N + d$
	    \EndWhile
	   \State  Return $N, C$ 
	\end{algorithmic} 
\end{algorithm} 

\begin{algorithm}\label{Distributed and Parallel computing Procedure}
	\caption{Distributed and Parallel Computing Procedure}
	\textbf{Initialize:} $M = 10$;  $N^c \leq 8$;  $P^t \leq 28$.
	\begin{algorithmic}[1]
	    \State \textbf{C} $\leftarrow$ 1
		\For {$j = 1,\ldots,M$}
		    \If{\textbf{C} $\leq N^c$} 
		        \State Assign $j$-th problem ($\hat{Q}$) to node \textbf{C}.
		        \State \textbf{C} $\leftarrow$ \textbf{C} + 1
		    \Else 
		        \State \textbf{C} $\leftarrow$ 1
		        \State Assign  $j$ problem ($\hat{Q}$) to node \textbf{C}.
		        \State \textbf{C} $\leftarrow$ \textbf{C} + 1
		    \EndIf
	    \EndFor
	    \For {$n = 1,\ldots,N^c$}
	        \If{Nr of problems assigned to n $>$ 1} 
	            \State Assign processors to problems
	            \State Solve ($\hat{Q}$) using Algorithm 1
	        \Else
	            \State Assign each of $P^t$ processors to the problem
	            \State Solve ($\hat{Q}$) using Algorithm 1
	        \EndIf
	    \EndFor
		\State Return  ${\bf x}_N$ for each of $M$ problems ($\hat{Q}$)
	\end{algorithmic} 
\end{algorithm} 

\clearpage

\onecolumn
\section{Input Data}
\begin{table*}[tbph!]
	\footnotesize
	\centering
	\begin{tabular}{lllccc}
		\toprule
		& & & \multicolumn{3}{c}{\textbf{Moisture Levels}} \\
		& & & \textbf{High} & \textbf{Medium} & \textbf{Low} \\
		\cline{4-6}
		\multirow{2}{*}{\textbf{Bale}} & \multicolumn{2}{l}{Moisture (\%)} & 25 & 17.50 & 10 \\
		& \multicolumn{2}{l}{Dry bulk density (dry Mg/cubic meter)} & 0.144 & 0.144 & 0.144 \\
		\midrule
		\multirow[b]{4}{*}{\textbf{First-stage}} & \multirow{3}{*}{Operating conditions} & Screen size (mm) & 76.2 & 76.2 & 76.2 \\
		& & Dry bulk density (dry Mg/cubic meter) & 0.144 & 0.144 & 0.144 \\
		& & Moisture (\%) & 25 & 17.5 & 10 \\
		\cline{2-6}
		& \multirow{5}{*}{Process performance} & Energy consumption (kWh/dry Mg) & 25.68 & 21.95 & 7.10 \\
		\multirow[t]{4}{*}{\textbf{Grinder}} & & Moisture loss (\%) & 4.77 & 3.00 & 0.5 \\
		& & Dry matter loss (\%) & 1.50 & 1.50 & 1.50 \\
		& & Bulk density change (dry Mg/cubic meter) & -0.09 & -0.10 & -0.11 \\
		& & Maximum in-feed rate (dry Mg/hour) & 2.20 & 4.53 & 5.23 \\
		\midrule
		\multirow[b]{3}{*}{\textbf{Separation}} & \multirow{3}{*}{Operating conditions} & Screen size (mm) & 6.35 & 6.35 & 6.35 \\
		& & Dry bulk density (dry Mg/cubic meter) & 0.053 & 0.041 & 0.039 \\
		& & Moisture (\%) & 20.23 & 14.50 & 9.50 \\
		\cline{2-6}
		\multirow[t]{2}{*}{\textbf{Unit}} & \multirow{2}{*}{Process performance} & Moisture loss (\%) & 0.71 & 0.50 & 0.00 \\
		& & Bypass ratio (\%) & 40.48 & 44.98 & 49.98 \\
		\midrule
		\multirow[b]{4}{*}{\textbf{Second-stage}} & \multirow{3}{*}{Operating conditions} & Screen size (mm) & 11.11 & 11.11 & 11.11 \\
		& & Dry bulk density (dry Mg/cubic meter) & 0.053 & 0.041 & 0.039 \\
		& & Moisture (\%) & 19.52 & 13.99 & 9.50 \\
		\cline{2-6}
		& \multirow{5}{*}{Process performance} & Energy consumption (kWh/dry Mg) & 49.71 & 16.58 & 14.67 \\
		\multirow[t]{4}{*}{\textbf{Grinder}} & & Moisture loss (\%) & 4.00 & 3.00 & 0.7 \\
		& & Dry matter loss (\%) & 0.50 & 0.50 & 0.50 \\
		& & Bulk density change (dry Mg/cubic meter) & 0.066 & 0.082 & 0.090 \\
		& & Maximum in-feed rate (dry Mg/hour) & 1.59 & 2.80 & 5.23 \\
		\midrule
		\multirow[b]{3}{*}{\textbf{Pelleting}} & \multirow{2}{*}{Operating conditions} & Dry bulk density (dry Mg/cubic meter) & 0.119 & 0.123 & 0.129 \\
		& & Moisture (\%) & 15.52 & 10.99 & 8.80 \\
		\cline{2-6}
		& \multirow{4}{*}{Process performance} & Energy consumption (kWh/dry Mg) & 90.39 & 60.63 & 55.12 \\
		\multirow[t]{3}{*}{\textbf{Equipment}} & & Moisture loss (\%) & 3.92 & 1.50 & 0.00 \\
		& & Bulk density change (dry Mg/cubic meter) & 0.547 & 0.542 & 0.537 \\
		& & Maximum in-feed rate (dry Mg/hour) & 3.33 & 3.81 & 4.76 \\
		\bottomrule
	\end{tabular}
	\caption{Technical process-related data for switchgrass}
	\label{tab:inputdata}
\end{table*}

\begin{table*}[tbph!]
  \footnotesize
  \centering 
  \caption{Short Sequence}
   \begin{tabular}{c|c|c}
\toprule      
{\bf Bale}  & \multicolumn{2}{c}{\bf Sequence Pattern}  \\  \cline{2-3} 
{\bf  Nr.}  & {\bf Feedstock-Based} & {\bf Moisture-Based} \\
\midrule
1-10   &        (2S-1C2)-(2S-3C2)-(1S-1C2)  &   3L-5M-2H \\
11-20  &        (1M-2C2)-(2S-3C2)-(1S-1C2)  &   3L-5M-2H \\
21-30    &      (1S-2C2)-(1M-1S-3C2)-(1M-1C2)   &   3L-5M-2H \\
31-40      &    (2C2-1C3)-(3C2-2C3)-(1C2-1C3)  &   3L-5M-2H \\
41-50   &       (1M-2C3)-(3C2-2C3)-(1C2-1C3)    &   3L-5M-2H \\
51-60    &      (1M-2C2)-(2M-1C2-2C3)-(2C2)     &   3L-5M-2H \\
61-70     &     (2C2-1C3)-(2M-3C3)-(1M-1C3)    &   3L-5M-2H \\
71-80      &    (1C2-2C3)-(4C2-1C3)-(1C2-1C3)   &   3L-5M-2H \\
    \hline
    \end{tabular}%
  \label{ShortSequence}%
\end{table*}%

\end{document}